\documentclass[12pt]{amsart}

\let\Box\qed
\let\fr\mathfrak
\let\frak\mathfrak
\let\cal\mathcal
\let\Bbb\mathbb
\def\Spin{\operatorname{Spin}}

\newtheorem{thm}{Theorem}[section] 
\newtheorem{prop}[thm]{Proposition} 
\newtheorem{cor}[thm]{Corollary} 
\newtheorem{lemma}[thm]{Lemma} 
\newtheorem{dof}{Definition}[section] 
\begin{document} 
\title
[$N$-(super)-extended Poincar\'e 
algebras and bilinear invariants]
{Classification of $N$-(super)-extended Poincar\'e 
algebras and bilinear invariants of the spinor representation of
$\Spin(p,q)$}
\author{Dmitry V. Alekseevsky}
\thanks{Alekseevsky was
supported  by the Max-Planck-Institut  f\"ur Mathematik (Bonn).
Cort\'es was
supported by the Alexander von Humboldt Foundation, MSRI (Berkeley)
and  SFB 256 (Bonn University). Research at MSRI is partially
supported through NSF grant DMS-9022140.}

\address{\hfuzz=5pt
Dmitry V. Alekseevsky\\
Max-Planck-Institut f\"ur Mathematik\\
Gottfried-Claren-Str.~26\\   
D-53225 Bonn}
\email{daleksee@@mpim-bonn.mpg.de}

\author{Vicente Cort\'es}
\address{Vicente Cort\'es\\
Mathematisches Institut der Universit\"at Bonn\\
Beringstra\ss e 1, D-53115 Bonn, Germany}

\curraddr{Mathematical Sciences
Research Institute\\
1000 Centennial Drive\\
Berkeley, CA 94720-5070} 
\email{vicente@@msri.org, vicente@@rhein.iam.uni-bonn.de} 

\begin{abstract} We classify extended Poincar\'e Lie  super algebras and 
Lie algebras of  any signature $(p, q)$, that is  
Lie super algebras  and ${\Bbb Z}_2$-graded Lie
algebras $\fr{g} =  \fr{g}_0 + \fr{g}_1$, where   
$\fr{g}_0 = {\fr{so}}(V) + V$ is the 
(generalized) Poincar\'e Lie algebra  of the pseudo Euclidean vector space  $ V = {\Bbb R}^{p,q}$ of   signature $(p,q)$
and  $\fr{g}_1 = S$ is the spinor ${\fr{so}}(V)$-module 
extended to a $\fr{g}_0$-module with kernel $V$.
The remaining super commutators $\{\fr{g}_1,\fr{g}_1\}$
(respectively, commutators $[\fr{g}_1, \fr{g}_1]$) 
are defined by an ${\fr {so}}(V)$-equivariant linear mapping
$$
\vee^2\fr{g}_1 \to V \quad (\mbox{respectively},\quad \wedge^2\fr{g}_1 \to V)\, .\hskip1in
$$
  Denote by ${\cal P}^+(n,s)$ (respectively, ${\cal P}^-(n,s)$) the vector 
space  of all such Lie super algebras (respectively, Lie algebras), where $n = p + q = \dim V$ and 
$s = p-q $ is the signature. The description of ${\cal P}^{\pm}(n,s)$ reduces 
to the construction of all ${\fr{so}}(V)$-invariant bilinear forms
on $S$  and to the calculation of three ${\Bbb Z}_2$-valued
invariants for some of them.

  This calculation is based on a simple explicit model of an irreducible
Clifford module  $S$  for the Clifford algebra $Cl_{p,q}$ of arbitrary
signature $(p,q)$.
As a result of the classification, we obtain the numbers
$L^{\pm}(n,s) = \dim {\cal P}^{\pm}(n,s)$
of independent  Lie super algebras and algebras, which 
take values 0,1,2,3,4 or 6.
Due to Bott periodicity, $L^{\pm}(n,s)$
may be considered as  periodic functions  with period 8  in
each argument. They are invariant under the group $\Gamma$
generated by  the four reflections with respect to the
axes $n=-2$, $n=2$, $s-1 = -2$ and $s-1 =2$. Moreover,
the reflection $(n,s) \to (-n,s)$ with respect to the
axis $n=0$  interchanges $L^+$ and $L^-$ :
$$  
L^+(-n,s) = L^-(n,s)\, .\hskip1in
$$ 
\end{abstract} 
\maketitle
\section*{Introduction}
General relativity is a gauge theory with the Poincar\'e
group $$P(1,3) = {\Bbb R}^{1,3} {\Bbb o} Lor(1,3)$$ of Minkowski 
space ${\Bbb R}^{1,3}$ as gauge group. In $N$-extended supergravity 
the $N$-extended Poincar\'e supergroup plays the role of (super) gauge 
group. 

The Lie super algebra of this super group for $N = 1$ is defined
as follows: $\fr{p}^{(1)}(1,3) = \fr{g} =  \fr{g}_0 + \fr{g}_1 = 
\fr{p}(1,3) + S$, where $\fr{p}(1,3) = {\Bbb R}^{1,3} + \fr{so}(1,3)$ 
is the Poincar\'e Lie algebra and $S = {\Bbb C}^2$ is the spinor module
of the Lorentz algebra $\fr{so}(1,3) \cong \fr{sl}(2,{\Bbb C})$ trivially extended to a 
$\fr{p}(1,3)$-module. The supercommutator ${\{ \cdot , \cdot \} } :S\otimes S 
\rightarrow {\Bbb R}^{1,3}$ is defined as projection onto the unique 
vector submodule $V \cong {\Bbb R}^{1,3}$ in the symmetric square $\vee^2S$.  

We remark that in this case there exists also a unique vector submodule in 
$\wedge^2S$, which defines on $\fr{p}(1,3) + S$ the structure of a 
${\Bbb Z}_2$-graded Lie algebra $\fr{p}^{(-1)}(1,3)$. 

Our goal is to classify for any pseudo Euclidean space $V = {\Bbb R}^{p,q}$ 
all similar extensions of the (generalized) Poincar\'e algebra $\fr{p}(V) = \fr{p}(p,q) 
= {\Bbb R}^{p,q} + \fr{so}(p,q)$ to a super Lie algebra or to a  
${\Bbb Z}_2$-graded Lie algebra. 

An other motivation to study such extensions is that extended Poincar\'e 
Lie algebras are closely related to the full isometry algebra $\fr{isom}(M)$ 
of homogeneous quaternionic K\"ahler manifolds $M$ (s.\ \cite{dW-V-VP}, \cite{A-C1}). 
In fact,  $\fr{isom}(M) = \fr{p} + {\Bbb R}A$, where $\fr{p}$ is an extension
of the Poincar\'e algebra $\fr{p}(3,3+k)$ of the pseudo Euclidean space 
${\Bbb R}^{3,3+k}$ of signature $(3,3+k)$, $k = -1,0,1, \ldots$, and $A$ is 
a derivation of $\fr p$ defining a natural gradation. 

\begin{dof}\label{N-extDef} A super Lie algebra (respectively a ${\Bbb Z}_2$-graded Lie algebra) 
$\fr{g} =  \fr{g}_0 + \fr{g}_1$ is called an {\bf $N$-extended} (respectively {\bf $-N$-extended}) 
{\bf Poincar\'e algebra} of $V = {\Bbb R}^{p,q}$ if the following conditions hold 
\begin{itemize}
\item[1)] $\fr{g}_0 \cong \fr{p}(V)$.  
\item[2)] $\fr{g}_1$ is a sum of $N$ irreducible spinor or semi spinor modules of 
$\fr{p}(V) = V + \fr{so}(V)$ with trivial action of the vector group V. 
\item[3)] The super bracket $\{ S,S\} \subset V$ (respectively Lie bracket $[S,S]\subset V$). 
\end{itemize}
\end{dof} 

Let $S$ be a $\fr{p}(V)$-module with trivial action of the vector group $V$. 
Then defining on $\fr{g} = \fr{p}(V) + S$ the structure of a super Lie algebra 
(respectively of a ${\Bbb Z}_2$-graded Lie algebra) such that 
$\fr{g}_0 \cong \fr{p}(V)$, $\fr{g}_1 = S$ and $\{ S,S\} \subset V$ (respectively
$[S,S]\subset V$) is equivalent to defining an $\fr{so}(V)$-equivariant 
mapping $j: V^{\ast} \rightarrow \vee^2S^{\ast}$ (respectively 
$j:V^{\ast} \rightarrow \wedge^2S^{\ast}$).  The super bracket (respectively
the Lie bracket) is given by $j^{\ast}:\vee^2S  \rightarrow V$ 
(respectively $j^{\ast}:\wedge^2S \rightarrow V$). Remark that under
these assumptions the Jacobi identities are automatically satisfied. 

We show that the classification of $N$-extended ($N\in {\Bbb Z}$) 
Poincar\'e algebras easily reduces to the classification of equivariant 
embeddings $V^{\ast} \hookrightarrow \vee^2S^{\ast}$ if $N>0$ and 
$V^{\ast} \hookrightarrow \wedge^2S^{\ast}$ if $N<0$, where $V$ is the 
vector module and $S$ the spinor module of $\fr{so}(V)$. In other 
words, we reduce the classification to the cases $N=\pm 1 ,\pm 2$. 

We prove that the following three vector spaces are isomorphic: 
\begin{itemize}
\item[1)] the space $\cal J$ of $\fr{so}(V)$-equivariant mappings 
$j: V^{\ast} \rightarrow  S^{\ast}\otimes S^{\ast}$,
\item[2)] the space $\cal M$ of  $\fr{so}(V)$-equivariant multiplications 
$\mu :  V^{\ast}\otimes S \rightarrow  S$ and 
\item[3)] the space $\cal B$ of $\fr{so}(V)$-invariant bilinear 
forms $\beta$ on $S$.    
\end{itemize} 

Let $\rho : V^{\ast}\otimes S\rightarrow  S$ be the (standard)
Clifford multiplication, where we have identified $V\cong V^{\ast}$ 
using the scalar product on $V = {\Bbb R}^{p,q}$. Then an 
isomorphism $j_{\rho}: {\cal B} \rightarrow  {\cal J}$
is given by 
\[ j_{\rho}(\beta ): v^{\ast}\in V^{\ast} \mapsto \beta \circ \rho (v^{\ast}) 
=\beta (\rho (v^{\ast})\cdot , \cdot )\in S^{\ast}\otimes S^{\ast}\, .\]

In particular, the classification of  $\fr{so}(V)$-equivariant mappings 
$V^{\ast} \rightarrow  S^{\ast}\otimes S^{\ast}$ is equivalent to the 
classification of $\fr{so}(V)$-invariant bilinear 
forms on the spinor module $S$. The latter amounts to the description 
of the Schur algebra $\cal C$ of $\fr{so}(V)$-invariant endomorphisms 
of $S$. The structure of $\cal C$ as abstract algebra depends only on the 
signature $s = p-q$ of ${\Bbb R}^{p,q}$ modulo 8; it is a simple 
real, complex or quaternionic matrix algebra of rank 1 or 2 or a sum of two 
isomorphic such algebras. 

To construct equivariant embeddings of the vector module $V^{\ast}$ 
into the symmetric square $\vee^2S^{\ast}$ (or into the exterior 
square $\wedge^2S^{\ast}$) we introduce the notion of  admissible
bilinear form $\beta$ on $S$ and also the corresponding notion of 
admissible endomorphism of $S$, which depends on the choice of an admissible
bilinear form $\beta$. 
\begin{dof} An $\fr{so}(V)$-invariant bilinear form $\beta$ on the spinor 
module $S$ is called 
{\bf admissible} if it   has the following properties:
\begin{itemize}
\item[1)] Clifford multiplication $\rho (v)$ is either $\beta$-symmetric
or $\beta$-skew symmetric. We define the type $\tau$ of $\beta$ to be $\tau (\beta )
= +1$ in the first case and $\tau (\beta )= -1$ in the second. 
\item[2)] $\beta$ is symmetric or skew symmetric. Accordingly, we 
define the symmetry $\sigma$ of $\beta$ to be  $\sigma (\beta ) = \pm 1$. 
\item[3)] If the spinor module is reducible, $S = S^{+} + S^{-}$, then
$S^{\pm}$ are either mutually orthogonal or isotropic. We put 
$\iota (\beta ) = +1$ in the first case, $\iota (\beta ) = -1$ in the 
second and call $\iota (\beta )$ the isotropy of $\beta$. 
\end{itemize}
\end{dof}    
Every admissible form $\beta$ defines an $\fr{so}(V)$-equivariant embedding 
$j_{\rho} (\beta ): V^{\ast} \rightarrow \vee^2 S^{\ast}$ if 
$\tau (\beta )\sigma (\beta) = +1$ or 
$j_{\rho} (\beta ): V^{\ast} \rightarrow \wedge^2 S^{\ast}$ if 
$\tau (\beta )\sigma (\beta) = -1$. Moreover, if $S = S^{+} + S^{-}$,
then either $S^{\pm}$ are orthogonal or isotropic for every bilinear 
form in the image of $j_{\rho}(\beta )$. 

The main part of the paper  
is the construction of an admissible basis for the space $\cal J$ of 
equivariant mappings $V^{\ast} \rightarrow S^{\ast}\otimes S^{\ast}$, 
i.e.\ a basis consisting of embeddings $j_{\rho} (\beta )$, where $\beta$ 
are admissible bilinear forms on $S$.

To describe all admissible forms $\beta$ we make use of very simple 
explicit models of the irreducible Clifford modules  inspired by Ra\u{s}evski\u{\i} \cite{R}.   
We prove that the problem reduces to the three fundamental cases $V = {\Bbb R}^{m,m}$, 
${\Bbb R}^{k,0}$ and ${\Bbb R}^{0,k}$ using the isomorphisms 
$C\! \ell_{m+k,m} \cong C\! \ell_{m,m} \hat{\otimes} C\! \ell_k$ and 
$C\! \ell_{m,m+k}\cong C\! \ell_{m,m}\hat{\otimes} C\! \ell_{0,k}$ and the 
algebraic properties of 
the fundamental invariants $\tau$, $\sigma$ and $\iota$ with respect to 
${{\Bbb Z}}_2$-graded tensor 
products. 

Moreover, we establish that for every pseudo Euclidean vector space 
$V = {\Bbb R}^{p,q}$ there is a prefered non degenerate 
$\fr{so}(V)$-invariant bilinear 
form $h$ on the spinor module $S$. This allows us to define {\em canonically}
the notion of admissible endomorphism of $S$ and the invariants 
$\tau$, $\sigma$ and $\iota$ for such endomorphisms.  They are 
multiplicative with respect to the composition 
$h\circ A = h (A \cdot, \cdot )$, 
$A \in {\cal C}$  admissible.

Finally, we explicitly construct in all the cases an admissible basis for the 
Schur algebra $\cal C$. This canonically yields admissible bases for the space
$\cal B$ of invariant bilinear forms and the space $\cal J$ of equivariant 
mappings. 

This gives an explicit description of all extended Poincar\'e algebras 
$\fr{g} = \fr{p}(V) +S$, where $S$ is the spinor module. The super
(respectively Lie) brackets $\vee^2S\rightarrow V$ (respectively 
$\wedge^2S\rightarrow V$) are given as linear combinations of mappings 
$j_i^{\ast}$, where the $j_i:V^{\ast} \rightarrow \vee^2S^{\ast}$ 
(respectively $V^{\ast} \rightarrow \wedge^2S^{\ast}$) form an admissible
basis for the space of $\fr{so}(V)$-equivariant mappings 
$V^{\ast} \rightarrow  \vee^2S^{\ast}$ (respectively 
$V^{\ast}\rightarrow  \wedge^2S^{\ast}$). 

If the spinor module $S$ is an irreducible $\fr{so}(V)$-module, we 
obtain all $N=\pm 1$ extended Poincar\'e algebras. If $S$ is reducible,
then we obtain  all $N=\pm 2$ extended Poincar\'e algebras and using 
the invariant $\iota$ we can determine all $N=\pm 1$ extended Poincar\'e algebras.
Sometimes there exist only trivial $N = 1$ (or $N = -1$) 
extended Poincar\'e algebras, i.e.\ $\{ S,S\} = 0$ (or $[S,S] = 0$). 

Given a pseudo Euclidean vector space $V = {\Bbb R}^{p,q}$, let $|N| = 1$ 
or $2$ denote the number of irreducible summands of the spinor module $S$ of 
$\fr{so}(V)$. For 
fixed $N = +|N|$ or $N = -|N|$ we give now        
the dimension $d_N$ of the vector space of  $N$-extended Poincar\'e 
algebra structures on  
$\fr{g} = \fr{p}(V) + S$.  

The function $d_N$, which  depends only on the signature $(p,q)$, admits  
a symmetry group $\Gamma$ generated by reflections. Moreover, there is an additional 
supersymmetry which relates the dimension $L^+ := d_{+|N|}$ of the space of super algebras   to the 
dimension $L^-:= d_{-|N|}$ of the space of Lie algebras. 
 
More precisely: Denote by $n = p+q$ the dimension and by $s = p-q$ the signature of 
$V = {\Bbb R}^{p,q}$ and let $L^+ = L^+(n,s)$ (respectively $L^-(n,s)$) be the
maximal number of linearly independent super algebra structures $\vee^2S \rightarrow 
V$ (respectively Lie algebra structures $\wedge^2S\rightarrow V$) on $\fr{g} = \fr{p}(V) + S$. 
The functions $L^+$ and $L^-$ are periodic with period 8 in each argument, hence we may 
consider them as functions on ${\Bbb Z}^2 =  {\Bbb Z}\times {\Bbb Z}$. 
The value of the pair $(L^+,L^-)$ is given in Table \ref{Tabelle}. 
\begin{table}[ht]\caption[Numbers $L^+$  of  super algebras and $L^-$ of  Lie algebras]{\label{Tabelle} 
The numbers $L^+$  of  super algebras and $L^-$ of  Lie algebras 
$\fr{g} = \fr{p}(V) + S$ are given as functions of the dimension $n$ and 
signature $s$ of $V$.  A fundamental domain for the reflection group $\Gamma$ is emphasized in bold face. The
supersymmetry axis is given by the equation $n = 0$.} 

\begin{centering}
\begin{tabular}{|r||c|c|c|c|c|c|c|c|c|}\hline 
$s$: & \multicolumn{9}{c|}{$(L^+(n,s),L^-(n,s))$}\\ \hline\hline 
5& &1,3& &1,3& &3,1& &3,1& \\ \hline
4&4,4& &2,6& &4,4& &6,2& &4,4\\ \hline
3& &1,3& &{\bf 1,3}& &{\bf 3,1}& & 3,1& \\ \hline
2&4,4& &{\bf 2,6}& &{\bf 4,4}& &{\bf 6,2}& &4,4\\ \hline
1&  &1,3& &{\bf 1,3}& &{\bf 3,1}& &3,1& \\ \hline
0&1,1& &{\bf 0,2}& &{\bf 1,1}&  &{\bf 2,0}& &1,1\\ \hline
$-1$& &0,1& &{\bf 0,1}& &{\bf 1,0}& &1,0& \\ \hline
$-2$&1,1& &0,2& &1,1&  &2,0& &1,1\\ \hline
$-3$&  &1,3& &1,3& &3,1& &3,1& \\ \hline\hline
$n$: &-4&-3&-2&-1&0&1&2&3&4\\ \hline   
\end{tabular}      

\end{centering}
\vskip18pt
\end{table}    
It follows from the inspection of this table, that the function $(L^+,L^-)$ is invariant under 
the group $\Gamma$ 
generated by the reflections with respect to the 4 axes  defined by the equations 
$n = -2$, $n = 2$, $s':= s-1 = -2$ and $s'= s-1 = 2$. A fundamental domain $F$ for $\Gamma$ 
is 
\[ F = \{ (n,s)\in {\Bbb Z}^2| -2\le n\le 2\, ,\quad
-2\le s' = s-1 \le 2\}\cap G\, , \] 
\[ G = \{ (n,s)|  \exists (p,q)\in {\Bbb Z}^2: n = p+q,
\quad s = p-q \} = \{ (n,s)\in {\Bbb Z}^2| n+s \quad \mbox{even}\} \]
and consists of 12 points. The values of the pair $(L^+,L^-)$ at these 
points are typed in bold face in Table \ref{Tabelle}.  

Moreover, the reflection $\theta$ with respect to  the axis $\{ n = 0\}$, 
$\theta :(n,s) \mapsto   (-n,s)$, is a supersymmetry of the pair $(L^+,L^-)$, 
that is it interchanges the number of Lie algebras and Lie super algebras: 
\[ (L^+(+n,s),L^-(+n,s)) = (L^-(-n,s),L^+(-n,s))\, .\]
Short:
\begin{displaymath} \fbox{$\displaystyle L^{\pm} = L^{\mp}\circ \theta$} \end{displaymath} 
A fundamental domain $\tilde{F}$ for the group $\tilde{\Gamma} = <\Gamma , \theta >$ is given by
\[ \tilde{F} = \{ (n,s) = (0,0),\: (0,2),\: (1,-1),\: (1,1),\: (1,3),\: (2,0),\: (2,2) \}\, . \] 
In terms of the coordinates $(p,q)$ a fundamental domain with $p\ge 0$ and  $q\ge 0$ is given by
\[ \tilde{D} =  \{ (p,q) =(2,0),\: (1,1),\: (3,0),\: (2,1),\: (1,2),\: (3,1),\: (2,2)\}\, . \]

\medskip\noindent 
{\bf Acknowledgements}\\ 
 The first author is very grateful to  Max-Planck-Institut f\"ur Mathematik 
for financial support and hospitality. 
The second author would like to thank S.-S.\ Chern and R.\ Osserman 
for inviting him to MSRI, where he is now enjoying his stay; he would also 
like to thank  W.\ Ballmann and U.\ Hamenst\"adt for encouragement and support.  
\newpage 
{\hfuzz=1in
\tableofcontents  }
\section[(Super) Extensions of the Poincar\'e algebra $\frak p(p,q)$ and 
$\Spin(p,q)$-equivariant embeddings ${{\Bbb R}}^{p,q} \hookrightarrow 
S^{\ast}\otimes  S^{\ast}$]{(Super) Extensions of the Poincar\'e
algebra $\frak p(p,q)$ and 
$\Spin(p,q)$-equi\-var\-iant embeddings ${{\Bbb R}}^{p,q} \rightarrow 
S^{\ast}\otimes  S^{\ast}$}
\subsection{Extending the  Poincar\'e algebra}
Let $V = {\Bbb R}^{p,q}$ be the pseudo Euclidean space with the metric
${<x,y>} = \sum_{i=1}^px^iy^i - \sum_{j=p+1}^{p+q}x^jy^j$.  We denote by
$\fr{so}(V) = \fr{so}(p,q)$ the pseudo ortho\-gonal Lie algebra and by  
$\fr{p}(V) = \fr{p}(p,q) = \fr{so}(V) + V$ the semidirect sum of  
$\fr{so}(V)$ and the Abelian ideal $V$, it is the Lie algebra of the isometry
group of $(V,{<\cdot ,\cdot >})$. We call $\fr{p}(V)$ the {\bf Poincar\'e 
algebra} of the space $V$. 
\begin{dof} A ${\Bbb Z}_2$-graded Lie algebra (respectively a super
algebra) $\fr{g} =  \fr{g}_0 + \fr{g}_1$ is called an {\bf extension} 
(respectively a {\bf super extension})
of $\fr{p}(V)$
 if $\fr{g}_0 = \fr{p}(V)$, 
$V$ is in the kernel of the representation of   $\fr{g}_0$ on 
$\fr{g}_1$ and $[\fr{g}_1,\fr{g}_1] \subset V$ (respectively 
$\{ \fr{g}_1,\fr{g}_1 \} \subset V$). 
\end{dof}
{\bf Remark 1:} Sometimes, for unification, we will refer to 
${\Bbb Z}_2$-graded Lie algebras and to super
algebras  as $\epsilon$-algebras, where $\epsilon = -1$ or 
$+1$ respectively. Correspondingly, we will speak of $\epsilon$-extensions. 

\begin{prop}\label{ext-pairProp} There exists a natural one-to-one correspondence between 
extensions (respectively super extensions) of $\fr{p}(V)$ up to isomorphisms 
and equivalence classes of pairs $(\rho , \pi )$, where 
\[ \rho : \fr{so}(V) \rightarrow \fr{gl}(W)\] 
is a representation and 
\[ \pi : \wedge^2W \rightarrow V \quad (\mbox{resp.} \quad \vee^2W  
\rightarrow V)\] 
is a $\fr{so}(V)$-equivariant linear map from  the space of skew symmetric
(respectively symmetric) bilinear forms on $W^*$ to the vector module $V$. 
Two pairs $(\rho , \pi )$ and $(\rho ', \pi ')$ ($\rho ' : \fr{so}(V) 
\rightarrow \fr{gl}(W')$) are {\bf equivalent} if there 
exists an automorphism $\phi :\fr{p}(V) \rightarrow \fr{p}(V)$ and a linear
map $\psi : W\rightarrow W'$ such that the following diagramms are commutative
(for  pairs of skew symmetric type): 
\[ \begin{array}{r@{\,}c@{\,}l@{\:}c@{\;}c@{\,}l} 
 &  \fr{so}(V)     &                                   &
\stackrel{\rho}{\longrightarrow} & \fr{gl}(V)          & \\ 
  & \downarrow &  {\scriptstyle \bar{\phi}}&                                    & \downarrow &  {\scriptstyle \psi}\\
&  \fr{so}(V)    &                                 &
\stackrel{\rho '}{\longrightarrow}   & \fr{gl}(W')          & 
\end{array} \qquad
\begin{array}{r@{\,}c@{\,}l@{\:}c@{\;}c@{\,}l} 
 &  \wedge^2W     &                                   &
\stackrel{\pi }{\longrightarrow} & V          & \\ 
  & \downarrow &  {\scriptstyle \psi}&                                    & 
\downarrow &  {\scriptstyle \phi|V}\\
&  \wedge^2W'    &                                &
\stackrel{\pi '}{\longrightarrow}   & V          & 
\end{array} \]
where $\bar{\phi}$ is the induced automorphism of $\fr{so}(V) = \fr{p}(V)/V$.
For pairs of symmetric type $\wedge^2$ must be replaced by $\vee^2$. 
\end{prop}
{\bf Proof:}
Given a pair $(\rho , \pi )$ of skew symmetric type, we define a   
${\Bbb Z}_2$-graded Lie algebra $\fr{g} =  \fr{g}_0 + {\fr 
g}_1$, $\fr{g}_0 = \fr{p}(V) = \fr{so}(V) + V$, ${\fr
g}_1 = W$ by 
\begin{eqnarray*}
{[}A,w{]} &=& \rho (A)w\, , \\
{[}w_1,w_2{]} &=& \pi (w_1\wedge w_2)\, ,\\ 
{[}v,w{]} &=& 0\, ,  
\end{eqnarray*} 
where $A\in \fr{so}(V)$, $v\in V$ and $w, w_1, w_2\in W$. For a pair 
of symmetric type we define a super algebra $\fr{g} =  \fr{g}_0 + {\fr 
g}_1$ by the same formulas replacing only the middle equation by 
\begin{eqnarray*} \{w_1,w_2\} &=& \pi (w_1\vee w_2)\, .
\end{eqnarray*} 
The Jacobi identity is satisfied because $\rho$ is a representation, 
$\pi$ is equivariant and the (anti)commutator of $W$ with $W$ is 
contained in $V$ and hence commutes with $W$. 
The other statements can be checked easily. $\Box$ 

Recall that the spinor representation is the representation  
of $\fr{so}(V)$ on an irreducible module $S$ of the Clifford 
algebra $C\! \ell (V)$. It is either irreducible or a sum of 
two irreducible semi spinor modules $S^{\pm}$. 

\begin{dof}\label{N-extDef2} (cf.\ Def.\ \ref{N-extDef}) Let $\fr{g} = \fr{g}(\rho , \pi )$ be an $\epsilon$-extension
of $\fr{p}(V)$ associated with a pair $(\rho , \pi )$. We say that $\fr g$ is 
an $\epsilon N$-{\bf extended Poincar\'e algebra} if $\rho$ is a sum of 
$N = 0,1,2,\ldots$ irreducible spin 1/2 representations, i.e.\ irreducible 
spinor or semi spinor representations. 
\end{dof} 

The purpose of this paper is to classify all  $N$-extended 
($N\in {\Bbb Z}$) Poincar\'e algebras. Before starting 
this classification we explain how, given a (super) extension of the 
Poincar\'e algebra, we can construct more complicated \linebreak[2] 
$\epsilon$-algebras. 
\subsection{Internal symmetries and  charges} 
\begin{dof}Let $\fr{g} = \fr{g}_0 + \fr{g}_1$ be an $\epsilon$-algebra.
An {\bf internal symmetry} of 
$\fr{g}$ is an automorphism of $\fr g$ which acts trivially on $\fr{g}_0$. 
\end{dof}
Now we give a simple construction which associates with an $\epsilon$-extension
$\fr{g} = \fr{g}(\rho ,\pi )$ of the Poincar\'e algebra $\fr{p}(V)$ and 
$l\in {\Bbb N}$ an $\epsilon$-extension  $\fr{g}^{(+l)}$ and also 
a $-\epsilon$-extension $\fr{g}^{(-2l)}$ which admit $O(l)$ respectively 
$Sp(2l,{\Bbb R})$ as internal symmetry groups.
We define $\fr{g}^{(+l)} = \fr{g}(\rho^{(+l)},\pi^{(+l)})$, where 
\[\rho^{(+l)} = l\rho: \fr{so}(V) \rightarrow lW = W\otimes {\Bbb R}^l
\,, \]
\[ \pi^{(+l)}(w_1\otimes v_1,w_2\otimes v_2)= \pi (w_1, w_2){<v_1,v_2>}\, ,\] 
$<\cdot ,\cdot >$ is the standard Euclidean scalar product on
${\Bbb R}^l$. Similarly, we define 
\[ \fr{g}^{(-2l)} = 2l\rho : \fr{so}(V) \rightarrow 2lW = W\otimes 
{\Bbb R}^{2l}\, ,\]
\[ \pi^{(-2l)}(w_1\otimes v_1,w_2\otimes v_2)= \pi (w_1, w_2)\omega (v_1,v_2) 
\, ,\]
where $\omega$ is the standard symplectic form on ${\Bbb R}^{2l}$. Here we have used
the convention that $\pi (w_1, w_2) = \pi (w_1\vee w_2)$ if $\epsilon = +1$ and 
$\pi (w_1, w_2) = \pi (w_1\wedge w_2)$ 
if  $\epsilon = -1$.  
\begin{prop} If $\fr{g}$ is an $\epsilon$-extension of the Poincar\'e 
algebra $\fr{p}(V)$, then $\fr{g}^{(+l)}$ is an $\epsilon$-extension
and $\fr{g}^{(-2l)}$ is a $-\epsilon$-extension. The standard actions of 
$O(l)$ (respectively $Sp(2l,{\Bbb R})$) on ${\Bbb R}^l$ 
(respectively ${\Bbb R}^{2l}$) are naturally extended to actions on 
$\fr{g}^{(+l)}$ (respectively $\fr{g}^{(-2l)}$) by internal symmetries. 
\end{prop} 
{\bf Proof:} 
The first statement follows immediately  from Prop.\ \ref{ext-pairProp} 
and the remark that the  bilinear map  $\pi^{(+l)}$ (respectively 
$\pi^{(-2l)}$) has the same (respectively the opposite) symmetry 
as $\pi$. The last statement is immediate. $\Box$ 

\noindent 
{\bf  Example 1:} Applying this construction to an $\epsilon$-extended 
(s.\ Def.\ \ref{N-extDef2}) Poincar\'e algebra, we obtain an 
$\epsilon l$-extended Poincar\'e algebra and also an  
$-\epsilon 2l$-extended Poincar\'e algebra with internal symmetry 
groups $O(l)$ and  $Sp(2l,{\Bbb R})$ respectively. 

\begin{dof} A ${\Bbb Z}_2$-graded Lie algebra (respectively 
a super algebra) $\fr{g} =  \fr{g}_0 + \fr{g}_1$ is called a 
{\bf charged extension} (respectively a {\bf charged super extension}) 
of the Poincar\'e algebra $\fr{p}(V)$ if 
\begin{itemize}
\item[1)]  $\fr{g}_0 = \fr{p}(V) + C$ is a trivial extension of 
$\fr{p}(V)$, i.e.\ $[C,C] = 0$. 
\item[2)] The action of $V + C$ on the 
$\fr{g}_0$-module $W = \fr{g}_1$ is trivial. 
\item[3)] The Lie (respectively super) bracket $\pi : \wedge^2W \rightarrow 
\fr{g}_0$ (respectively $\vee^2W \rightarrow 
\fr{g}_0$) is a sum $\pi = \pi_V + \pi_C$, where 
$\pi_V: \wedge^2W \rightarrow V$ and $\pi_C: \wedge^2W \rightarrow C$ 
(respectively $\pi_V: \vee^2W  \rightarrow V$ and $\pi_C: \vee^2W \rightarrow 
C$). In particular, $(\fr{p}(V) + W,\pi_V)$ is  an extension (respectively 
super extension) of $\fr{p}(V)$. 
\end{itemize} 
If moreover, $[\fr{so}(V),C] = 0$, and hence $[C,\fr{g}] = 0$, then 
$\fr{g}$ is called a {\bf central charge  extension} (respectively 
a {\bf central charge  super extension}) of  $\fr{p}(V)$. 
\end{dof} 

Let an extension (respectively super extension)  $\fr{p}(V) + W$ admitting 
a connected Lie group $H$ of internal symmetries be given.  Without 
restriction of generality we can assume that 
$H$ is simply connected and we denote the Lie algebra of 
$H$ by $\fr h$. To construct a charged  extension 
(respectively super extension) $(\fr{p}(V) + C) + W$ preserving the 
internal symmetry group $H$ it is necessary and 
sufficient to define an $(\fr{so}(V)+\fr{h})$-equivariant map
$\pi_C$ from 
the exterior (respectively symmetric) square of $W$ to an  
$(\fr{so}(V)+\fr{h})$-module $C$. 

\noindent
{\bf Example 2:} Let $\fr{p}(V) + W$ be an extension of $\fr{p}(V)$. 
Consider the extension $\fr{g}^{(+l)} = \fr{p}(V) + W
\otimes {\Bbb R}^l$ with internal symmetry group $H = O(l)$ defined 
above. Let $h\in \vee^2 W^{\ast}\otimes {\Bbb R}^r$  be a  
symmetric $\fr{so}(V)$-invariant 
(possibly trivial) vector valued bilinear form on $W$ and $\eta  
\in \wedge^2 W^{\ast}\otimes {\Bbb R}^s$   a  skew 
symmetric such form. Define 
\[ \pi_C: \wedge^2(W\otimes {\Bbb R}^l) \rightarrow C = 
{\Bbb R}^r\otimes \wedge^2{\Bbb R}^l + {\Bbb R}^s\otimes 
 \vee^2{\Bbb R}^l\, ,\] 
\[\pi_C(w_1\otimes x_1,w_2\otimes x_2) = h(w_1,w_2)x_1\wedge x_2 + 
\eta (w_1,w_2)x_1\vee x_2\, ,\] 
where $w_1, w_2\in W$ and $x_1, x_2 \in  {\Bbb R}^l$. Then  $\pi_C$ 
defines on $(\fr{p}(V) + C) + W\otimes {\Bbb R}^l$ the structure 
of central charge extension of $\fr{p}(V)$  with symmetry group $O(l)$. 

Analogeously, we can define on 
$(\fr{p}(V) + C) + W\otimes {\Bbb R}^{2l}$, $C = {\Bbb R}^r 
\otimes \vee^2{\Bbb R}^{2l} +  {\Bbb R}^s \otimes 
\wedge^2{\Bbb R}^{2l}$, 
the structure of central charge super extension of $\fr{p}(V)$  with 
symmetry group $Sp(2l,{\Bbb R})$ by 
\[ \pi_C: \vee^2(W\otimes {\Bbb R}^{2l}) \rightarrow C \, ,\]    
\[\pi_C(w_1\otimes x_1,w_2\otimes x_2) = h(w_1,w_2)x_1\vee x_2 + 
\eta (w_1,w_2)x_1\wedge x_2\, .\]  

\noindent
{\bf Example 3:} Let $\fr{p}(V) + W$ be a super extension of $\fr{p}(V)$. 
Consider the super extension $\fr{g}^{(+l)} = \fr{p}(V) + 
W\otimes {\Bbb R}^{l}$ 
with internal symmetry group $H = O(l)$ and let $h$ be a symmetric and $\eta$ 
a skew symmetric  vector valued 
$\fr{so}(V)$-invariant  bilinear form on $W$, as above. Define 
\[ \pi_C: \vee^2(W\otimes {\Bbb R}^{l}) \rightarrow C = {\Bbb R}^r 
\otimes \vee^2{\Bbb R}^{l} + {\Bbb R}^s \otimes 
\wedge^2{\Bbb R}^{l}\, ,\] 
\[\pi_C(w_1\otimes x_1,w_2\otimes x_2) = h(w_1,w_2)x_1\vee x_2 + 
\eta (w_1,w_2)x_1\wedge x_2\, .\] 
Then  $\pi_C$ 
defines on $(\fr{p}(V) + C) + W\otimes {\Bbb R}^l$ the structure 
of central charge super extension of $\fr{p}(V)$  with symmetry group $O(l)$.
 
Analogeously, we can define on 
$(\fr{p}(V) + C) + W\otimes {\Bbb R}^{2l}$, $C = {\Bbb R}^r 
\otimes \wedge^2{\Bbb R}^{2l} + 
{\Bbb R}^s \otimes \vee^2{\Bbb R}^{2l}$ 
the structure of central charge extension of $\fr{p}(V)$  with 
symmetry group $Sp(2l,{\Bbb R})$ by 
\[ \pi_C: \wedge^2(W\otimes {\Bbb R}^{2l}) \rightarrow C \, ,\] 
\[\pi_C(w_1\otimes x_1,w_2\otimes x_2) = h(w_1,w_2)x_1\wedge x_2 + 
\eta (w_1,w_2)x_1\vee x_2\, .\]

In the physical literature (s.\ \cite{F}) the expression ``central 
charges'' is used for a special case of Example 3.  
\subsection{Reduction of the classification of $N$-extended Poin\-ca\-r\'e
algebras to the cases $N = \pm 1, \pm 2$} 
Let $\fr{g} = \fr{g}(\rho ,\pi ) = \fr{p}(V) + W$  be a $\pm N$-extended 
Poincar\'e algebra, $N = 1,2,\ldots$. Then either the spinor representation 
$\rho_0 : \fr{so}(V)\rightarrow \fr{gl}(S)$ is irreducible and 
$\rho = N\rho_0$, $W = NS = S\otimes {\Bbb R}^N$, or it decomposes
into two irreducible subrepresentations $\rho_0 = \rho_+ + \rho_-$, 
$S = S^+ +  S^-$ and $\rho = N_+\rho_+ + N_-\rho_-$, 
$W = N_+S^+ + N_-S^- = S^+\otimes {\Bbb R}^{N_+} + 
S^-\otimes {\Bbb R}^{N_-}$, $N = N_+ + N_-$. 
The description of all $\epsilon N$-extended Poincar\'e algebras 
$\fr{g}(\rho ,\pi )$ reduces to the description of all 
$\fr{so}(V)$-equivariant mappings $\pi :\wedge^2W\rightarrow V$ 
if $\epsilon = -1$ and   $\pi :\vee^2W\rightarrow V$ if $\epsilon = +1$. 
If $\pi \neq 0$, the dual mapping defines an  $\fr{so}(V)$-equivariant 
embedding $\pi^{\ast}: V^{\ast} \hookrightarrow \wedge^2W^{\ast}$ 
if $\epsilon = -1$ or   $\pi^{\ast}: V^{\ast} \hookrightarrow \vee^2W^{\ast}$ 
if $\epsilon = +1$. To find all such embeddings it is suficient to determine
all submodules isomorphic to $V^{\ast}$ in $\wedge^2W^{\ast}$ and 
$\vee^2W^{\ast}$ or, equivalently, all vector submodules  $V$
in $\wedge^2W$ and $\vee^2W$. Tables \ref{vee2WTab} and \ref{wedge2WTab} 
reduce this problem to the 
cases $N = 1$ or $2$.  
\begin{table}[ht]\caption[Decomposition of the   symmetric square 
of $W$]{\label{vee2WTab} Decomposition of the   symmetric square 
of $W$}

\begin{centering}
\begin{tabular}{|c||c|c|}\hline 
$\rho$:  & $N\rho_0$ & $N_+\rho_+ + N_-\rho_-$  \\ \hline 
$W$:     & $NS = S\otimes {\Bbb R}^N$ & $N_+S^+ + N_-S^- =$   \\ 
 & & $S^+\otimes 
{\Bbb R}^{N_+} + S^-\otimes {\Bbb R}^{N_-}$ \\ \hline 
$\vee^2W$ & 
$\vee^2S\otimes \vee^2{\Bbb R}^N 
+  \wedge^2S\otimes \wedge^2{\Bbb R}^N$  &  
$\vee^2S^+\otimes \vee^2{\Bbb R}^{N_+} + \vee^2S^-\otimes 
\vee^2{\Bbb R}^{N_-} +$ 
   \\ 
 &  & 
 $\wedge^2S^+\otimes 
\wedge^2{\Bbb R}^{N_+} + \wedge^2S^-\otimes 
\wedge^2{\Bbb R}^{N_-} +$ 
  \\
 & & $S^+\otimes S^-\otimes {\Bbb R}^{N_+N_-}$ \\  
\hline 
\end{tabular} 

\end{centering}
\vskip18pt
\end{table}
\begin{table}[ht]\caption[Decomposition of the  exterior square 
of $W$]{\label{wedge2WTab} Decomposition of the   exterior square 
of $W$}

\begin{centering}
\begin{tabular}{|c||c|c|}\hline 
$\rho$:  & $N\rho_0$ & $N_+\rho_+ + N_-\rho_-$  \\ \hline 
$W$:     & $NS = S\otimes {\Bbb R}^N$ & $N_+S^+ + N_-S^- =$   \\ 
 & & $S^+\otimes 
{\Bbb R}^{N_+} + S^-\otimes {\Bbb R}^{N_-}$ \\ \hline 
  
$\wedge^2W$&      
$\wedge^2S\otimes \vee^2{\Bbb R}^N +  \vee^2S\otimes \wedge^2{\Bbb R}^N$   
 & $\wedge^2S^+\otimes \vee^2{\Bbb R}^{N_+} + \wedge^2S^-\otimes 
\vee^2{\Bbb R}^{N_-} + $\\ 

 &  & $ \vee^2S^+\otimes 
\wedge^2{\Bbb R}^{N_+} + \vee^2S^-\otimes 
\wedge^2{\Bbb R}^{N_-} + $\\
 & & $S^+\otimes S^-\otimes {\Bbb R}^{N_+N_-}$\\ 
\hline
\end{tabular} 

\end{centering}
\vskip18pt
\end{table}

If $\rho_+$ and $\rho_-$ are equivalent then $\rho = N_+\rho_+ 
+N_-\rho_- \cong N\rho_0$, $\rho_0 \cong \rho_{\pm}$, 
\begin{eqnarray*} \vee^2W 
& \cong & \vee^2S_0\otimes \vee^2{\Bbb R}^N + 
\wedge^2S_0\otimes \wedge^2{\Bbb R}^N\, ,\\
\wedge^2W 
& \cong & \vee^2S_0\otimes \wedge^2{\Bbb R}^N + 
\wedge^2S_0\otimes \vee^2{\Bbb R}^N\, ,
\end{eqnarray*} 
where $S_0 \cong S^{\pm}$ and $N = N_+ +N_-$. Table \ref{vee2WTab}  
shows that the classification of all equivariant embeddings 
$V\hookrightarrow \vee^2W$ (case $\epsilon = +1$) reduces to finding 
all equivariant embeddings $V\hookrightarrow \vee^2S$ and 
$V\hookrightarrow \wedge^2S$ if $S$ is irreducible and 
equivariant embeddings $V\hookrightarrow \vee^2S^{\pm}$, 
$V\hookrightarrow \wedge^2S^{\pm}$ and  $V\hookrightarrow S^+\otimes S^-$ 
if $S = S^+ + S^-$. Table \ref{wedge2WTab} shows that the same reduction 
applies to  the case $\epsilon = -1$, i.e.\ to the problem of finding all 
equivariant embeddings $V\hookrightarrow \wedge^2S$. We see that e.g. 
the classification of $N$-extended Poincar\'e algebras for $N>0$ (i.e.\  
super algebra extensions) reduces to the classification of 
$N = \pm 1$-extended Poincar\'e algebras in  case there is only one 
irreducible spin 1/2 representation of $\fr{so}(V)$. The same is true 
for $N<0$, i.e.\ for Lie algebra extensions. 

To illustrate this reduction we consider the case $\epsilon = +1$ and 
$\rho = N\rho_0$ in more detail. 
\begin{lemma} Assume $\epsilon = +1$ and 
$\rho = N\rho_0$, where $\rho_0$ is an irreducible spin 1/2 representation 
on $S_0$. Then any $\fr{so}(V)$-equivariant embedding 
\[ j: V\hookrightarrow \vee^2W = \vee^2S_0\otimes \vee^2{\Bbb R}^N + 
\wedge^2S_0\otimes \wedge^2{\Bbb R}^N\] 
is given by 
\[ j(v) = \sum_a\phi_a(v)\otimes A_a + \sum_b \psi_b(v)\otimes B_b\, ,\] 
where $\phi_a: V\rightarrow \vee^2S_0$ and $\psi_b: V\rightarrow \wedge^2S_0$ 
are equivariant embeddings, $A_a\in \vee^2{\Bbb R}^N$ and 
$B_b\in \wedge^2{\Bbb R}^N$.  
\end{lemma}   
{\bf Proof:} Choose bases $(A_a)$ and $(B_b)$ of $\vee^2{\Bbb R}^N$ 
and $\wedge^2{\Bbb R}^N$ respectively. Then $j(v)$ can be decomposed as 
above and the coefficients $\phi_a$ and $\psi_b$ are equi\-variant embeddings 
or zero.  $\Box$ 

\subsection{Equivariant embeddings $V^{\ast}\hookrightarrow S^{\ast} \otimes 
S^{\ast}$ , modified Clifford multiplications and  Dirac operators} 
We reduced the problem of the classification of $N$-extended Poincar\'e 
algebras to the description of $\fr{so}(V)$-equivariant mappings 
$V^{\ast} \rightarrow S^{\ast} \otimes S^{\ast}$, where $S$ is 
the spinor  module of $\fr{so}(V)$. We will denote by $\cal J$
the vector space of all such mappings. 

Now we will show that this space is closely related to two other vector spaces:
\begin{itemize}
\item the space $\cal B$ of all $\fr{so}(V)$-invariant bilinear forms 
on $S$ and 
\item the space $\cal M$ of $\fr{so}(V)$-equivariant multiplications $\mu : 
V^{\ast}\otimes S\rightarrow S$. 
\end{itemize} 
Denote by $\cal C$ the {\bf Schur algebra} of $\fr{so}(V)$-invariant 
endomorphisms of $S$. We define two natural anti-representations of 
$\cal C$ on $\cal B$ and $\cal J$ and also a representation and an 
anti-representation of $\cal C$ on $\cal M$ by: 
\begin{eqnarray*}
\xi_A^{\cal B}\beta & = & \beta (A\cdot , \cdot )\\
\eta_A^{\cal B}\beta & = & \beta (\cdot , A\cdot )\\
(\xi_A^{\cal J}j)(v^{\ast}) & = & \xi_A^{\cal B}(j(v^{\ast}))\\
(\eta_A^{\cal J}j)(v^{\ast}) & = & \eta_A^{\cal B}(j(v^{\ast}))\\
(\xi_A^{\cal M}\mu )(v^{\ast}) & = & A\circ \mu (v^{\ast})\\ 
(\eta_A^{\cal M}\mu )(v^{\ast}) & = & \mu (v^{\ast})\circ A\, ,
\end{eqnarray*} 
where $A\in {\cal C}$, $v^{\ast}\in V^{\ast}$, $\beta \in {\cal B}$, 
$j\in {\cal J}$ and $\mu \in {\cal M} \subset Hom(V^{\ast}, End\, S)$. 
Remark that a non zero equivariant mapping $j: V^{\ast}\rightarrow 
S^{\ast}\otimes S^{\ast}$ is automatically an embedding. 

\begin{dof} An equivariant embedding $j: V^{\ast}\rightarrow 
S^{\ast}\otimes S^{\ast}$ is called {\bf non degenerate}, if 
$j(V^{\ast})S = S^{\ast}$ and $j(S) \cong S$, where we consider 
$j$ as mapping $j: S\rightarrow V\otimes S^{\ast}$.  An 
equivariant multiplication $\mu : V^{\ast}\otimes S\rightarrow  S$ is called 
{\bf non degenerate}, if $\mu (V^{\ast})S = S$. 
\end{dof} 

Using the following identifications, we define mappings from two of the spaces ${\cal B}$, $\cal J$ and $\cal M$ into the third. 
\begin{eqnarray*} {\cal B} & = &  
(S^{\ast}\otimes S^{\ast})^{\fr{so}(V)}\, ,\\ 
 {\cal J} & = & Hom(V^{\ast},S^{\ast}\otimes S^{\ast})^{\fr{so}(V)} 
\stackrel{(\ast )}{\cong}  
Hom (S,V^{\ast}\otimes S^{\ast})^{\fr{so}(V)}\, ,\\
 {\cal M} & = &  Hom(V^{\ast}\otimes S,S)^{\fr{so}(V)}
 \cong 
 Hom(V^{\ast},End\, S)^{\fr{so}(V)}\\ & \cong & Hom(V^{\ast}
\otimes S^{\ast},S^{\ast})^{\fr{so}(V)}\, .
\end{eqnarray*}
at ($\ast$)  we used the metric identification $V^{\ast} \cong V$.
The mappings are defined as follows: 
\begin{eqnarray*}{\cal B}\times {\cal M} & \rightarrow & {\cal J} \\
  (\beta ,\mu ) & \mapsto & j(\beta ,\mu ) = \beta \circ \mu \\
j(\beta ,\mu )(v^{\ast}) &=& \beta (\mu (v^{\ast}) \cdot , \cdot )\, , \quad 
v^{\ast}\in V^{\ast}\, ;\\ 
{\cal M}\times {\cal J} & \rightarrow & {\cal B}\\ 
(\mu ,j) & \mapsto & \beta(\mu , j) = \mu \circ j\, ,\\
\beta(\mu , j)(s,t) &=& {< \mu (j(s)),t>}\, , \quad s,t \in S\, ;\\ 
{\cal B}\times {\cal J} & \rightarrow & {\cal M}\\ 
(\beta ,j) & \mapsto & \mu (\beta ,j) = \beta \circ j\\ 
\mu (\beta ,j)(v^{\ast}) &=& \beta (j(v^{\ast}) \cdot , \cdot )\in S\otimes 
S^{\ast} \cong End\, S\, , 
\end{eqnarray*}
where $<\cdot ,\cdot >$ denotes the natural duality pairing $S^{\ast}\times 
S\rightarrow {\Bbb R}$ and for the last mapping we have used that 
$j(v^{\ast}) \in S^{\ast}\otimes S^{\ast} 
\cong Hom (S^{\ast},S)$. 

\begin{thm}\label{BJMT} The choice of a non degenerate element $\beta_0$, 
$j_0$ or $\mu_0$ in any of the spaces 
${\cal B}$, $\cal J$ and $\cal M$ defines vector space isomorphisms between 
the two  others: 
\begin{eqnarray*} j_{\beta_0}: {\cal M} & \rightarrow & {\cal J}\\
 \mu     & \mapsto & j(\beta_0,\mu ) = \beta_0 \circ \mu\, , \\
\mu_{\beta_0}:  {\cal J} & \rightarrow &  {\cal M}\\ 
j & \mapsto &\mu (\beta_0,j) = \beta_0 \circ j\, ;\\ 
\beta_{j_0}: {\cal M} & \rightarrow &  {\cal B}\\
\mu &  \mapsto &\beta (\mu , j_0) = \mu \circ j_0\, ,\\
\mu_{j_0}: {\cal B} & \rightarrow & {\cal M}\\
\beta & \mapsto & \mu (\beta ,j_0) = \beta \circ j_0\, ;\\
j_{\mu_0} : {\cal B} & \rightarrow &  {\cal J}\\ 
\beta & \mapsto & j(\beta , \mu_0) = \beta \circ \mu_0\, ,\\
\beta_{\mu_0} : {\cal J} & \rightarrow & {\cal B}\\ 
j & \mapsto & \beta (\mu_0 , j) = \mu_0 \circ j\, .
\end{eqnarray*} 
\end{thm} 
{\bf Proof:} The statement is trivial for $j_{\beta_0}$ and $\mu_{\beta_0}$,
because these mappings amount to ``raising and lowering'' indices of 
tensors via the non degenerate form $\beta_0$. 

It is clear that $\mu_{j_0}$ and $j_{\mu_0}$ are injective, since $j_0$ and 
$\mu_0$ are non degenerate. Hence, it is sufficient to prove that 
$\beta_{j_0}$ and $\beta_{\mu_0}$ are injective. 

Consider first $\beta_{\mu_0}(j) = \mu_0 \circ j$, where $j:S\rightarrow 
V^{\ast}\otimes S^{\ast}$ and $\mu_0: V^{\ast}\otimes S^{\ast}\rightarrow 
S^{\ast}$. The kernel of $\beta_{\mu_0}$
equals 
\[ ker \, \beta_{\mu_0} = \{j\in{\cal J}| j(S)\subset ker\, \mu_0\} \]  
If $0\neq j\in ker \, \beta_{\mu_0}$, then $ker\, \mu_0$ contains the 
non trivial submodule $j(S)$. This is impossible, 
because  $ker\, \mu_0$ does not contain spin 1/2 submodules. Indeed, 
after complexification the $\fr{so}(V^{{\Bbb C}})$-module 
$(V^{\ast})^{{\Bbb C}}\otimes (S^{\ast})^{{\Bbb C}}$ has the 
decomposition 
\[ (V^{\ast})^{{\Bbb C}}\otimes (S^{\ast})^{{\Bbb C}} = \Sigma 
\oplus (S^{\ast})^{{\Bbb C}} = (ker\, \mu_0^{{\Bbb C}})\oplus 
(S^{\ast})^{{\Bbb C}}\, ,\] 
where $\Sigma = ker\, \mu_0^{{\Bbb C}}$ contains only spin 3/2 
modules, i.e.\ Kronecker product of the vector module $V^{{\Bbb C}}
\cong (V^{\ast})^{{\Bbb C}}$ 
(spin 1) and an irreducible spin 1/2 module. 

Consider now $\beta_{j_0}(\mu ) = \mu \circ j_0$, where $j_0: S\rightarrow 
V^{\ast}\otimes S^{\ast}$ and $\mu : V^{\ast}\otimes S^{\ast}\rightarrow 
S^{\ast}$. As before we have the decomposition 
$ (V^{\ast})^{{\Bbb C}}\otimes (S^{\ast})^{{\Bbb C}} = \Sigma 
\oplus (S^{\ast})^{{\Bbb C}}$, where $\Sigma$ has no submodules 
isomorphic to submodules of $(S^{\ast})^{{\Bbb C}}$. 
If $\mu \neq 0$, $ker\, \mu^{{\Bbb C}} = \Sigma \oplus 
S_1^{{\Bbb C}}$, where $S_1^{{\Bbb C}} \neq 
(S^{\ast})^{{\Bbb C}}$ is a proper submodule of 
$(S^{\ast})^{{\Bbb C}}$. Since $j_0$ is non degenerate 
$j_0(S) \cong S$ cannot be contained in $ker\, \mu$. $\Box$ 

\begin{lemma}\label{selfdualL}  Let $S$ be the spinor module of $\fr{so}(V)$.
There always exists a non degenerate $\fr{so}(V)$-invariant bilinear form 
$\beta$ on $S$. 
\end{lemma}
{\bf Proof:} The existence of $\beta$ is equivalent to the self duality 
of $S$, i.e.\ to the condition $S^{\ast} \cong S$ as $\fr{so}(V)$-modules. 

The self duality    of the complex $\fr{so}(V^{{\Bbb C}})$ spinor 
module ${{\Bbb S}}$ follows from the criterion of self duality given in 
\cite{O-V} p.\ 195. 

Now we discuss the real case. Assume first  $S^{{\Bbb C}}$ has the same 
number of irreducible summands as $S$. Then the self duality of $S$ follows 
from that of $S^{{\Bbb C}}$, s.\ \cite{O-V} p.\ 291. In the opposite case
$S$ admits an invariant complex structure $J$ and $(S,J) \cong {\Bbb S}$
(complex spinor module of $\fr{so}(V^{{\Bbb C}})$). 
Then the real part of a non degenerate complex 
$\fr{so}(V^{{\Bbb C}})$-invariant bilinear form on $S = {\Bbb S}$ 
gives a real $\fr{so}(V)$-invariant bilinear form on $S$ and it is easy to 
check that this form is non degenerate. $\Box$ 

From Theorem \ref{BJMT} and this lemma we now derive an important consequence. 
Recall that by definition the spinor module $S$ is a module over the Clifford 
algebra $C\!\ell (V)$. The restriction of the multiplication mapping 
$C\!\ell (V) \times S \rightarrow S$ to $V\times S$ defines a non degenerate   
$\fr{so}(V)$-equivariant multiplication $\rho : V\otimes S \cong V^{\ast}
\otimes S \rightarrow S$, 
which is called Clifford multiplication (as above $V$ and $V^{\ast}$ are 
identified using the pseudo Euclidean scalar product of $V$). 
The composition $j(\beta ,\rho ) =
\beta \circ \rho$ with a non degenerate $\fr{so}(V)$-invariant form $\beta$ 
gives a non degenerate $\fr{so}(V)$-equivariant embedding 
$V^{\ast} \hookrightarrow S^{\ast}\otimes S^{\ast}$.  
Using the lemma and this remark, we obtain the following corollary from 
Theorem \ref{BJMT}.  

\begin{cor}\label{BIMCor} The spaces ${\cal B}$ of $\fr{so}(V)$-invariant bilinear forms on 
$S$, $\cal J$ of $\fr{so}(V)$-equivariant mappings 
$V^{\ast} \rightarrow S^{\ast}\otimes S^{\ast}$ and $\cal M$ of 
$\fr{so}(V)$-equivariant mul\-ti\-pli\-ca\-tions 
$V^{\ast} \otimes S \rightarrow S$ are isomorphic. In particular, 
Clifford multiplication $\rho$ defines the isomorphism 
$j_{\rho}: {\cal B} \rightarrow {\cal J}$ and hence any 
$\fr{so}(V)$-equivariant embedding 
$V^{\ast} \hookrightarrow S^{\ast}\otimes S^{\ast}$ is of the form 
\[ j = j_{\rho}(\beta ): v^{\ast}\mapsto \beta (\rho (v^{\ast})\cdot , \cdot )
\, , \quad \beta \in {\cal B}\, , \quad v^{\ast}\in V^{\ast}\, .\]
\end{cor} 
{\bf Remark 2:} Using an $\fr{so}(V)$-equivariant multiplications
$\mu :V^{\ast}\otimes S \rightarrow S$ one can define a 
Dirac type operator $D^{\mu}$ on a pseudo Riemannian spin manifold $M$ as 
follows. 
Let $\mu_x: T^{\ast}_xM\otimes S_x \rightarrow S_x$ be a field of equivariant 
multiplications, where $S(M) = \cup_{x\in M}S_x \rightarrow M$ is the 
spinor bundle. Then 
\[ (D^{\mu}s)_x = \mu_x( \nabla s) = \mu_x (\sum_ie^i\otimes \nabla_{e_i}s)\, ,\] 
where $(e_i)$ is a basis of $T_xM$, $(e^i)$ the dual basis of $T^{\ast}_xM$ 
and $\nabla$ is the spinor connection  induced by the Levi Civita connection. 
\subsection{${{\Bbb Z}}_2$-graded type and Schur algebra $\cal C$} 
It is well known (s.\ \cite{L-M}), that every Clifford algebra 
$C\!\ell (V)$, $V = 
{\Bbb R}^{p,q}$, is isomorphic to ${\Bbb K}(l)$ or to 
$2{\Bbb K}(l) = {\Bbb K}(l)\oplus {\Bbb K}(l)$, where 
${\Bbb K}(l)$ is the full matrix algebra over ${\Bbb K}$ of
rank $l$ depending on $(p,q)$ and   where ${\Bbb K} = {\Bbb R}$, ${{\Bbb C}}$ or ${{\Bbb H}}$. 
\begin{dof} \label{typeDef} We say that a Clifford algebra $C\!\ell (V)$ 
has {\bf type} 
$r{\Bbb K}$, $r = 1$ or $2$, if $C\!\ell (V) \cong r {\Bbb K}(l)$ 
for some $l\in {\Bbb N}$. 
\end{dof} 
  
Recall that the Clifford  algebra $C\!\ell (V)$ has a natural 
${{\Bbb Z}}_2$-grading $C\!\ell (V) = C\!\ell^0 (V) + C\!\ell^1 (V)$. If 
$V = {\Bbb R}^{p,q}$ ($\neq 0$), then the even part  $C\!\ell^0 (V)$ 
is isomorphic 
to the Clifford algebra $C\!\ell (V')$ of $V' = {\Bbb R}^{p-1,q}$ if 
$p\ge 1$ and $V' = {\Bbb R}^{q-1}$ if $p = 0$. 
Remark that $\dim C\!\ell^0 (V)= {\dim C\!\ell (V)}/2$. By  
the preceding remarks, the following definition makes sense. 
\begin{dof}\label{Z2typeDef} The pair 
$$t(C\!\ell (V)) = (r_0{\Bbb K}_0,r{\Bbb K}) =
({type} \, C\!\ell^0 (V), {type} \, C\!\ell (V))$$ is called the 
$\Bbb Z_2$-{\bf graded type} of the Clifford  algebra $C\!\ell (V)$. 
\end{dof}

The following proposition describes the periodicity of the type $t$ of the 
$\Bbb Z_2$-graded  Clifford  algebras $C\!\ell_{p,q} = 
C\!\ell ({\Bbb R}^{p,q})$. 
\begin{prop}\label{typeProp} The $\Bbb Z_2$-graded type 
$t_{p,q} = t(C\!\ell_{p,q})$ depends only on the si\-gna\-ture $s = p-q$ 
modulo 8 and $t(s) = t(p-q) = t_{p,q}$ is given in the 
table.  \\ 

\begin{centering}
\begin{tabular}{|c||c|c|c|c|c|c|c|c|}\hline 
$s$ & $1$ & $2$ & $3$& $4$& $5$ & $6$ & $7$ & $8$\\ \hline\hline 
$t(s)$ & ${\Bbb R},{\Bbb C}$&${\Bbb C},{\Bbb H}$&
${\Bbb H},2{\Bbb H}$& $2{\Bbb H},{\Bbb H}$ & 
${\Bbb H},{\Bbb C}$ & ${\Bbb C},{\Bbb R}$& 
${\Bbb R},2{\Bbb R}$& $2{\Bbb R},{\Bbb R}$\\ \hline 
\end{tabular} 

\end{centering}
\end{prop}
{\bf Proof:} The proof reduces to the investigation of \cite{L-M} Table II. 
$\Box$ 

\begin{cor}The $\Bbb Z_2$-graded type $t_{p,q} = t(s=p-q)$ is mirror 
symmetric with respect to the diagonal $\{ p + q = 0\}$: 
$t_{p,q} = t_{-q,-p}$; in other words,  $t(C\!\ell_{p,q}) = 
t(C\!\ell_{8k-q,8k-p})$, $8k\ge p,q$. 

Moreover, the $\Bbb Z_2$-graded type $t_{p,q} = t(s) = 
(t^0(s),t^1(s))$ is mirror super symmetric with respect to the 
axis 
$\{ s = p-q = 3.5\}$, i.e.\ 
\[ (t^0(7-s),t^1(7-s)) = (t^1(s),t^0(s))\, .\]     
\end{cor}
The type $r{\Bbb C}$ and $\Bbb Z_2$-graded type 
$t_m = (r_0{\Bbb C},r{\Bbb C})$ of a complex Clifford algebra 
$C\!\ell_m = C\!\ell ({\Bbb C}^m)$ are defined by putting $V = 
{\Bbb C}^m$ in Definition \ref{typeDef} and \ref{Z2typeDef}, where  
${\Bbb C}^m$ is equipped with a non degenerate (complex) bilinear form, 
e.g.\ the standard one: 
$<z,w> = \sum_{j=1}^m z_jw_j$, $z,w \in {\Bbb C}^m$. 
\begin{prop}
The $\Bbb Z_2$-graded type 
$t_m = t({\Bbb C}\!\ell_m)$ depends only on the parity of $m$: 
\[ t_m = \left\{ \begin{array}{r@{\quad \mbox{if} \quad}l}
(2{\Bbb C},{\Bbb C})& \mbox{{\it m} is even}\\
({\Bbb C},2{\Bbb C})&  \mbox{{\it m} is odd} 
\end{array}      \right.\] 
\end{prop}   

Let $S = S_{p,q}$ be an irreducible  $C\!\ell_{p,q}$-module. Recall that 
by definition the Schur algebra ${\cal C} = {\cal C}_{p,q}$ of $S$ is the 
algebra of all its   $\fr{so}(V)$-invariant endomorphisms; it is the algebra 
of endomorphisms which commute with $C\!\ell^0_{p,q}$. Analogously, we define 
the Schur algebra ${\cal C}^c_m$ of the complex spinor module ${{\Bbb S}}$; it is 
the algebra of endomorphism of ${\Bbb S}$ commuting with 
${\Bbb C}\!\ell^0_m$. 

\begin{cor}\label{SchurCor} The Schur algebra ${\cal C}_{p,q} = {\cal C}(p-q)$ depends only 
on $s = p-q$ modulo 8 and is given in the table. In particular, 
it admits the mirror symmetry $(p,q) \mapsto (-q,-p)$.\\  

\begin{centering}
\begin{tabular}{|c||c|c|c|c|c|c|c|c|}\hline 
$s$ & $1$ & $2$ & $3$& $4$& $5$ & $6$ & $7$ & $8$\\ \hline\hline  
${\cal C}(s)$ & ${\Bbb R}(2)$ & ${\Bbb C}(2)$ & 
${\Bbb H}$ & ${\Bbb H} \oplus {\Bbb H}$ & ${\Bbb H}$  & 
${\Bbb C}$ & ${\Bbb R}$ & ${\Bbb R} \oplus {\Bbb R}$  
\\ \hline 
\end{tabular} 

\end{centering}
\end{cor} 
{\bf Proof:} Remark that if $t(C\!\ell_{p,q}) = 
(r_0 {\Bbb K}_0,r{\Bbb K})$ and hence $C\!\ell^0_{p,q} \cong 
r_0 \, {\Bbb K}_0(l_0)$, $C\!\ell_{p,q} \cong 
r{\Bbb K}(l)$, then $l$ is completely determined by $l_0$ and vice versa; 
$l = l_0$ or $2l_0$.  This follows from $\dim C\!\ell_{p,q} = 
2\dim C\!\ell^0_{p,q}$. 

Using this remark, Proposition \ref{typeProp} 
shows that the pair $(Cl^0_{p,q}, Cl_{p,q})$ is isomorphic to one of the following: 
\begin{eqnarray*} 
({\Bbb K}(l),{\Bbb K}'(l)) &,& \quad S = {{\Bbb K}'}^l\, ,\\ 
({\Bbb K}(l),2{\Bbb K}(l)) &,& \quad S = {\Bbb K}^l\, ,\\
({\Bbb K}'(l),{\Bbb K}(2l))&,& \quad S = {\Bbb K}^{2l}\, ,\\
(2{\Bbb K}(l), {\Bbb K}(2l))&,& \quad S = {\Bbb K}^{2l}\, ,
\end{eqnarray*} 
where ${\Bbb K} = {\Bbb R}$, ${\Bbb C}$ or ${\Bbb H}$ 
and ${\Bbb R}' = {\Bbb C}$, ${\Bbb C}' = {\Bbb H}$. 

In the first case the ${\Bbb K}(l)$-module $S = {{\Bbb K}'}^l$ 
is a sum of two irreducible equivalent modules 
$S^{\pm} \cong  {\Bbb K}^l$ and hence the Schur algebra  ${\cal C} 
\cong {\Bbb K}(2)$.  

In the second (respectively third) case $S =  {\Bbb K}^l$ 
(respectively ${\Bbb K}^{2l}$) is irreducible as 
${\Bbb K}(l)$- (respectively ${\Bbb K}'(l)$-) module and hence 
${\cal C} \cong {\Bbb K}$ (respectively ${\Bbb K}'$). 

In the last case ${\cal C} \cong {\Bbb K}\oplus {\Bbb K}$, 
which follows from the next lemma.    

\begin{lemma} Let $S = {\Bbb K}^{2l}$ be the irreducible module of the 
algebra  ${\Bbb K}(2l)$ and ${\cal A} \cong 2{\Bbb K}(l)$ a 
subalgebra of ${\Bbb K}(2l)$, then the ${\cal A}$-module $S$ is 
decomposed  into a sum of two nonequivalent submodules $S^{\pm}$.
\end{lemma}
{\bf Proof:} It is clear that the ${\cal A}$-module $S$ is the sum of two 
irreducible submodules $S^+$ and $S^-$. They are not equivalent because 
${\cal A}|S^{+}$ and ${\cal A}|S^{-}$ have different kernels, 
namely the two ideals ${\Bbb K}(l)\subset {\cal A}$. $\Box$

Remark that the algebras ${\Bbb C}\oplus {\Bbb C}$ and 
${\Bbb H}(2)$ do not occur as Schur algebras of the real spinor 
module $S$. 

\begin{cor} The Schur algebra ${\cal C}^c_m$ of the complex spinor module 
${{\Bbb S}}$ depends only 
on the parity of $m$:
\[ {\cal C}^c_m \: = \: \left\{ \begin{array}{r@{\quad \mbox{if} \quad}l}
{\Bbb C}\oplus {\Bbb C}& \mbox{{\it m} is even}\\
{\Bbb C}&  \mbox{{\it m} is odd} 
\end{array}      \right.\]  
\end{cor} 

The proof of Corollary \ref{SchurCor} shows that the structure of the 
matrix algebra  ${\cal C}$ contains the following information about the 
$C\!\ell^0(V)$-module $S$. 
\begin{prop} ${\cal C}$ is a simple ${\Bbb K}$-matrix algebra (respectively a 
sum of two isomorphic ${\Bbb K}$-matrix algebras) if and only if $C\!\ell^0(V)$ 
is a simple ${\Bbb K}$-matrix algebra (respectively a 
sum of two isomorphic such algebras). $S$ is an irreducible 
$C\!\ell^0(V)$-module  if and only if ${\cal C} \cong  
{\Bbb K}$ ($= {\Bbb R}$, ${\Bbb C}$ or ${\Bbb H}$). 
$S$ is decomposed into a sum 
of two equivalent (respectively inequivalent) $C\!\ell^0(V)$-modules  if and 
only if ${\cal C} \cong {\Bbb K}(2)$ (respectively 
${\cal C} \cong {\Bbb K}\oplus {\Bbb K}$). 
\end{prop} 

The corresponding statement in the complex case is given for the sake of 
completeness: 
\begin{prop}If $m$ is even, then the spinor module ${\Bbb S} = 
{\Bbb S}_m$ is the sum ${\Bbb S} = {\Bbb S}^+ + 
{\Bbb S}^-$ of two inequivalent irreducible 
${\Bbb C}\!\ell^0_m$-modules. In this case, 
${\Bbb C}\!\ell^0_m$
and the Schur algebra ${\cal C}^c_m$  are the direct sum  of two isomorphic  
simple (complex) matrix algebras. 

If $m$ is odd, then the spinor module is an irreducible module of the 
simple  matrix algebra ${\Bbb C}\!\ell^0_m$ and its Schur algebra
is also simple.
\end{prop}       
 
Since, due to  Lemma \ref{selfdualL}, $S$ admits a non degenerate 
$\fr{so}(p,q)$-invariant bilinear form, by Schur's Lemma the dimension
$b_{p,q}$ of the space ${\cal B} = {\cal B}_{p,q}$ of  
$\fr{so}(p,q)$-invariant bilinear forms on $S$ equals 
\[ b_{p,q} = \dim {\cal B}_{p,q} = \dim {\cal C}_{p,q}\, .\] 
Hence we have: 
\begin{cor}  \label{bpqCor} $b_{p,q} = b(p-q)$ is a periodic function of 
$s = p-q$ with period 8. In particular, it admits the mirror 
symmetry   $(p,q) \mapsto (-q,-p)$. Its values are given in the following 
table.\\ 

\begin{centering}
\begin{tabular}{|c||c|c|c|c|c|c|c|c|}\hline 
$s$ & $1$ & $2$ & $3$& $4$& $5$ & $6$ & $7$ & $8$\\ \hline\hline  
$b(s)$ & 4 & 8 & 4 & 8 &4  & 2  & 1 & 2  
\\ \hline 
\end{tabular} 

\end{centering}
\end{cor}   

Denote by $b_m$ the (complex) dimension of the space of 
$\fr{so}(m,{\Bbb C})$-invariant bilinear forms on the complex 
spinor module ${\Bbb S}$, then $b_m = \dim_{{\Bbb C}} {\cal C}^c_m$ and 
we have:
\[ b_m \: = \:  \left\{ \begin{array}{r@{\quad \mbox{if} \quad}l}
 2 & \mbox{{\it m} is even}\\
 1 &  \mbox{{\it m} is odd.} 
\end{array}      \right. \]

\section{Fundamental invariants $\tau$ , $\sigma$ and $\iota$ and  reduction 
to the  basic signatures $(m,m)$, $(k,0)$ and $(0,k)$}
\subsection{Fundamental invariants} 
As before let $V$ denote a pseudo Euclidean vector space and $S$ its spinor 
module.  In Corollary \ref{BIMCor} we have established that every 
$\fr{so}(V)$-equivariant embedding 
$j: V^{\ast} \hookrightarrow S^{\ast}\otimes 
S^{\ast}$ is of the form
\[ j = j_{\rho}(\beta ): v^{\ast}\mapsto \beta (\rho (v^{\ast})\cdot , \cdot )
\, ,\quad v^{\ast}\in V^{\ast}\, ,\]
where $\rho$ is Clifford multiplication and  $\beta \in {\cal B}$. 
The dimension of the space $\cal B$ of $\fr{so}(V)$-invariant bilinear forms
on $S$ was given in Corollary \ref{bpqCor}. 

Now we will concentrate on a class of bilinear forms $\beta \in {\cal B}$
for which  $j_{\rho}(\beta )V^{\ast}\subset \vee^2S^{\ast}$ or 
$j_{\rho}(\beta )V^{\ast}\subset \wedge^2S^{\ast}$ and define fundamental 
invariants $\tau$, $\sigma$ and $\iota$ for this class. 
\begin{dof} A bilinear form $\beta$ on the spinor module $S$ is called 
{\bf admissible} if it has the following properties:
\begin{itemize}
\item[1)] Clifford multiplication $\rho (v)$, $v\in V$, is either $\beta$-symmetric
or $\beta$-skew symmetric. We define the {\bf type} $\tau$ of $\beta$ to be $\tau (\beta )
= +1$ in the first case and $\tau (\beta )= -1$ in the second. 
\item[2)] The bilinear form $\beta$ is symmetric or skew symmetric. 
Accordingly, we 
define the {\bf symmetry} $\sigma$ of $\beta$ to be  
$\sigma (\beta ) = \pm 1$. 
\item[3)] If the spinor module is reducible, $S = S^{+} + S^{-}$, then
$S^{\pm}$ are either mutually orthogonal or isotropic. We put 
$\iota (\beta ) = +1$ in the first case, $\iota (\beta ) = -1$ in the 
second and call $\iota (\beta )$ the {\bf isotropy} of $\beta$. 
\end{itemize}
\end{dof} 
Due to 1) every admissible form $\beta$ is $\fr{so}(V)$-invariant and 
hence defines an $\fr{so}(V)$-equivariant embedding 
$j_{\rho}(\beta ): V\cong V^{\ast} \hookrightarrow S^{\ast}\otimes S^{\ast}$. 
In addition, 
$j_{\rho}(\beta )V\subset \vee^2S^{\ast}$ if 
$\tau (\beta )\sigma (\beta ) = +1$ and 
$j_{\rho}(\beta )V\subset \wedge^2S^{\ast}$ if 
$\tau (\beta )\sigma (\beta ) = -1$. If $S = S^+ + S^-$, 
then for every bilinear form $\gamma \in j_{\rho}(\beta )V$ 
the semi spinor modules $S^{\pm}$ are either $\gamma$-isotropic
(if $\iota (\gamma ) = -\iota (\beta ) = -1$) or mutually  
$\gamma$-orthogonal (if $\iota (\gamma ) = -\iota (\beta ) = +1$). 

Given an admissible form $\beta \in {\cal B}$ and $A\in {\cal C}$ the 
composition $\beta \circ A = \beta (A\cdot , \cdot ) \in {\cal B}$ is 
in general not admissible. However, if $A$ is $\beta$-admissible (see
Definition \ref{betaadmDef} below) then $\beta \circ A$ is admissible. 
\begin{dof}\label{betaadmDef} Let $\beta \in {\cal B}$ be admissible. An endomorphism 
$A$ of $S$ is called $\beta$-{\bf admissible} if it has the following 
properties:
\begin{itemize}
\item[1)] Clifford multiplication $\rho (v)$, $v\in V$, either commutes or 
anticommutes with $A$. We  define the {\bf type} $\tau$ of $A$ to be 
$\tau (A)
= +1$ in the first case and $\tau (A)= -1$ in the second. 
\item[2)] $A$ is $\beta$-symmetric or $\beta$-skew symmetric. Accordingly, we 
define the $\beta$-{\bf symmetry} $\sigma$ of $A$ to be  
$\sigma_{\beta} (A) = \pm 1$.{\hfuzz=2pt\par}
\item[3)] If the spinor module is reducible, $S = S^{+} + S^{-}$, then
either $AS^{\pm} \subset S^{\pm}$ or  $AS^{\pm} \subset 
S^{\mp}$. We put 
$\iota (A) = +1$ in the first case, $\iota (A) = -1$ in the 
second and call $\iota (A)$ the {\bf isotropy} of $A$.
\end{itemize}
\end{dof}
Due to 1) every $\beta$-admissible endomorphism $A$ is $\fr{so}(V)$-invariant
and hence $\beta \circ A\in {\cal B}$. Moreover, $\beta \circ A$ is 
admissible and the fundamental invariants are multiplicative: 
\begin{eqnarray*} 
\tau (\beta \circ A) &=&  \tau (\beta )\tau (A)\, ,\\
\sigma (\beta \circ A) &=&  \sigma (\beta )\sigma (A)\, ,\\
\iota (\beta \circ A) &=&  \iota (\beta )\iota (A)\, .
\end{eqnarray*}

In section \ref{signmmsubSec} (s.\ Definition \ref{canDef}), for every pseudo Euclidean
space $V$, we will construct a canonical non degenerate  
$\fr{so}(V)$-invariant  bilinear form $h$ on the spinor module $S$. 
We will define that an endomorphism $A$ of $S$ is admissible of symmetry 
$\sigma (A) = \pm 1$, if $A$ is $h$-admissible and $\sigma_h(A) = \pm 1$. 

\noindent 
{\bf Remark 3:} The complete classification of admissible forms $\beta \in 
{\cal B}$, which we will give in this paper, implies the following. 
Let $\gamma \in {\cal B}$ be  non degenerate and admissible.
Then a $\gamma$-admissible endomorphism $A\in {\cal C}$  is 
$\beta$-admissible for every admissible $\beta \in {\cal B}$.  
In particular, admissibility (i.e.\ $h$-admissibility) implies 
$\beta$-admissibility. 

\subsection{Reduction to the basic signatures}\label{redbasicsignSec} 
Let $V_1$ and $V_2$ be pseudo Euclidean spaces and $V = V_1 + V_2$ 
their orthogonal sum. We recall (s.\ \cite{L-M} I. Prop.\ 1.5) that 
there is a canonical isomorphism of $\Bbb Z_2$-graded algebras 
\[ C\!\ell (V) \cong C\!\ell (V_1) \hat{\otimes} C\!\ell (V_2) \, ,\]
where $\hat{\otimes}$ denotes the $\Bbb Z_2$-graded tensor product of 
$\Bbb Z_2$-graded algebras. 
\begin{prop} \label{Z2redProp} Let $M_1 = M_1^0 + M_1^1$ be a $\Bbb Z_2$-graded 
$C\!\ell (V_1)$-module and $M_2$ a (not necessarily $\Bbb Z_2$-graded)  
$C\!\ell (V_2)$-module. 
Then $M = M_1 \otimes M_2$ carries a natural structure of 
$C\!\ell (V)$-module, $V = V_1 +V_2$, given by: 
\[ (a_1\otimes a_2)(m_1\otimes m_2) = (-1)^{\deg (a_2) \deg (m_1)} 
{a_1m_1} \otimes {a_2m_2}\, ,\] 
where $a_i\in C\!\ell (V_i)$, $m_i\in M_i$, $i=1,2$.  
If $M_2 = M_2^0 + M_2^1$ is a  $\Bbb Z_2$-graded $C\!\ell (V_2)$-module,
then this formula defines on $M$ the structure of  $\Bbb Z_2$-graded 
$C\!\ell (V)$-module:  $M^0 = M_1^0\otimes M_2^0 + M_1^1\otimes M_2^1$, 
$M^1 = M_1^0\otimes M_2^1 + M_1^1\otimes M_2^0$.  
\end{prop}
\begin{cor}\label{redCor} Let $S_i$ be an irreducible $C\!\ell (V_i)$-module, $i=1,2$, 
and assume that $S_1 = S_1^+ + S_1^-$ is reducible as 
$C\!\ell^0 (V_1)$-module. Then $S = S_1\otimes S_2$ is an irreducible 
($C\!\ell (V) = C\!\ell (V_1) \hat{\otimes} C\!\ell (V_2)$)-module. 
The $C\!\ell^0(V)$-module $S$ is reducible, $S = S^+ + S^-$, if 
and only if $S_2$ is reducible as $C\!\ell^0(V_2)$-module, $S_2 = 
S_2^+ + S_2^-$. 
\end{cor} 

\noindent
{\bf Proof:} Let $S_1$ be an irreducible $C\!\ell (V_1)$-module 
which is reducible as $C\!\ell^0 (V_1)$-module and let $S_1^+$ be an 
irreducible $C\!\ell^0 (V_1)$-submodule. Then 
\[ S_1':= C\!\ell (V_1)\otimes_{C\!\ell^0 (V_1)}S_1^+\] 
is an irreducible $C\!\ell (V_1)$-module, hence $S_1 \cong S_1'$ as 
$C\!\ell (V_1)$-modules. Moreover, $S_1'$ is a  $\Bbb Z_2$-graded 
$C\!\ell (V_1)$-module (s.\ \cite{L-M} I.\ Prop.\ 5.20): $S_1' = {S_1'}^0 
+ {S_1'}^1$, ${S_1'}^0 = C\!\ell^0(V_1)\otimes_{C\!\ell^0 (V_1)}S_1^+ \cong 
S_1^+$ and ${S_1'}^1 = C\!\ell^1(V_1){S_1'}^0 = C\!\ell^1(V_1)
\otimes_{C\!\ell^0 (V_1)}S_1^+$. 

Therefore, we may assume (as usual) that 
$S_1 = S_1^+ + S_1^-$ is a 
$\Bbb Z_2$-graded $C\!\ell (V_1)$-module: $S_1^0 = S_1^+$, $S_1^1 = 
S_1^-= C\!\ell^1 (V_1)S_1^+$, reducing the first statement to Proposition
\ref{Z2redProp}. The remaining statements also follow from the structure of 
$\Bbb Z_2$-graded Clifford module on $S_1$ and on $S_2$ (in the reducible 
case).  $\Box$ 

Now we investigate the algebraic properties of the fundamental invariants 
with respect to $\Bbb Z_2$-graded tensor products. 
\begin{prop}\label{multiplProp} Under the assumptions of Corollary \ref{redCor} 
let $\beta_i$ be admissible bilinear forms on $S_i$, $i=1,2$. 

If $\tau (\beta_1) = \iota (\beta_1)\tau (\beta_2)$, then 
$\beta = \beta_1 \otimes \beta_2$ is admissible and 
\begin{eqnarray*}
\tau (\beta ) &=& \tau (\beta_1) \: =\:  \iota (\beta_1)\tau (\beta_2)\, ,\\
\sigma (\beta ) &=& \sigma (\beta_1)\sigma (\beta_2)\, ,\\
\iota (\beta ) &=& \iota (\beta_1)\iota (\beta_2)\, ,
\end{eqnarray*}
where $\iota (\beta )$ and $\iota (\beta_2)$ are defined if and only if 
$S_2$ (and hence $S$) is reducible as module of the even part of the
corresponding Clifford algebra. 

Let $A_i$ be $\beta_i$-admissible endomorphisms of $S_i$, $i=1,2$. 
If $\tau (A_1) = \iota (A_1)\tau (A_2)$, then $A = A_1 \otimes A_2$ 
is admissible and 
\begin{eqnarray*}
\tau (A) &=&  \tau (A_1) \: =\: \iota (A_1)\tau (A_2)\, ,\\
\sigma_{\beta}(A) &=& \sigma_{\beta_1}(A_1)\sigma_{\beta_2}(A_2)\, ,\\ 
\iota (A) &=& \iota (A_1) \iota (A_2)\, ,
\end{eqnarray*}
where $\iota (A)$ and $\iota (A_2)$ are defined if and only if 
$S_2$ is reducible as $C\!\ell^0(V_2)$-module.
\end{prop}

\noindent
{\bf Proof:} The only non trivial statements are the ones concerning the 
type $\tau$. For $s_i$, $t_i\in S_i$ and $v_i\in V_i$ we compute: 
\begin{eqnarray*} 
\beta ((v_1 \otimes 1)(s_1 \otimes s_2), t_1 \otimes t_2) &=& 
\beta (v_1s_1 \otimes s_2, t_1\otimes t_2)\: =\\
\beta_1(v_1s_1,t_1)\beta_2( s_2,t_2) &=& \tau (\beta_1)\beta_1(s_1,v_1t_1)
\beta_2( s_2,t_2)\: =\\
\tau (\beta_1)\beta (s_1 \otimes s_2,v_1t_1\otimes t_2) &=& 
\tau (\beta_1)\beta (s_1 \otimes s_2,(v_1 \otimes 1)(t_1 \otimes t_2))
\end{eqnarray*} 
and 
\begin{eqnarray*} 
\beta ((1\otimes v_2)(s_1 \otimes s_2), t_1 \otimes t_2) &=& 
(-1)^{\deg s_1} \beta (s_1 \otimes v_2s_2, t_1 \otimes t_2)\: =\\
(-1)^{\deg s_1} \beta_1(s_1,t_1)\beta_2(v_2s_2,t_2)  &=& 
(-1)^{\deg s_1}\tau (\beta_2)\beta_1(s_1,t_1)\beta_2(s_2,v_2t_2)\: =
\end{eqnarray*} 
\[ (-1)^{\deg s_1}\tau (\beta_2)\beta (s_1 \otimes s_2,t_1 \otimes v_2t_2)
\: =\] 
\[ (-1)^{\deg s_1 + \deg t_1}\tau (\beta_2)\beta (s_1 \otimes s_2,
(1\otimes v_2)( t_1 \otimes t_2))\, .\]

If $\iota (\beta_1) = (-1)^{\deg s_1 + \deg t_1}$ we obtain 
\begin{equation}\label{betaEqu} 
\beta ((1\otimes v_2)(s_1 \otimes s_2), t_1 \otimes t_2) = 
\iota (\beta_1)\tau (\beta_2)\beta (s_1 \otimes s_2,
(1\otimes v_2)( t_1 \otimes t_2))\, . 
\end{equation}
Otherwise, both sides of  (\ref{betaEqu}) vanish. Hence, the equation
(\ref{betaEqu}) is always true. 

Similarly we have: 
 \[ (v_1 \otimes 1)((A_1 \otimes A_2) (s_1\otimes s_2)) =
\tau (A_1) (A_1 \otimes A_2)((v_1 \otimes 1)(s_1\otimes s_2))\] 
and 
\begin{eqnarray*} 
(1\otimes v_2)((A_1 \otimes A_2) (s_1\otimes s_2)) &=& 
(1\otimes v_2)(A_1s_1\otimes  A_2s_2)\: =\\ 
(-1)^{\deg (A_1s_1)} A_1s_1\otimes v_2A_2s_2 &=& 
(-1)^{\deg (A_1s_1)}\tau (A_2) A_1s_1\otimes A_2v_2s_2\: =
\end{eqnarray*}  
\[ (-1)^{\deg (A_1s_1)}\tau (A_2)(A_1 \otimes A_2)(s_1\otimes v_2s_2) \: =\] 
\[ (-1)^{\deg (A_1s_1) + \deg s_1}\tau (A_2)(A_1 \otimes A_2)((1\otimes v_2)
(s_1\otimes s_2))\: =\]  
\[ \iota (A_1) \tau (A_2) (A_1 \otimes A_2)((1\otimes v_2)
(s_1\otimes s_2))\, . \quad \Box\]

Now we point out that every pseudo Euclidean space $V$ can be decomposed 
as orthogonal sum $V = V_1 + V_2$ such that the assumptions of 
Corollary \ref{redCor} are satisfied, i.e.\ such that the spinor 
$C\!\ell^0 (V_1)$-module $S_1$ is reducible. In fact, we can 
decompose $V$ into $V_1 = {\Bbb R}^{m,m}$ and $V_2 =  {\Bbb R}^{k,0}$
or ${\Bbb R}^{0,k}$. 
\begin{prop}\label{redProp} Let $V = V_1 + V_2$ be the orthogonal sum of the pseudo 
Euclidean spaces $V_1 = {\Bbb R}^{m,m}$ and $V_2$. Let $S_1$ be an 
irreducible $C\!\ell (V_1)$-module. Then $S_1 = S_1^+ + S_1^-$ is a sum of two
inequivalent irreducible $C\!\ell^0 (V_1)$-submodules $S_1^{\pm}$ and  
an irreducible 
($C\!\ell (V) = C\!\ell (V_1)\hat{\otimes}C\!\ell (V_2)$)-module 
$S$ is given by $S = S_1 \otimes S_2$, where $S_2$ is an irreducible 
$C\!\ell (V_2)$-module. $S$ is reducible as  $C\!\ell^0(V)$-module 
if and only if $S_2$ is reducible as $C\!\ell^0(V_2)$-module. 
\end{prop}

\noindent
{\bf Proof:}  The first statement follows from the fact that the Schur 
algebra of $S_1$ is ${\cal C}_{m,m} = {\cal C} (s = m-m = 0) = {\Bbb R}
\oplus {\Bbb R}$. Now all other statements follow immediately from 
Corollary \ref{redCor}. $\Box$

\section{Case of signature $(m,m)$ and complex case} 
\subsection{Signature $(m,m)$}\label{signmmsubSec}
Let $U$ and $U^{\ast}$ denote two complementary isotropic subspaces of 
$V = {\Bbb R}^{m,m}$, so 
$V = U + U^{\ast}$. We denote by ${<\cdot ,\cdot >}$ the scalar product of 
$V$ and identify $U^{\ast}$ with the dual space to $U$ 
by 
\[ u^{\ast}(u) \: =\:  {<u,u^{\ast}>}\, , \quad u^{\ast}\in U^{\ast}\, , \: 
u\in U\, .\]
\begin{prop} \label{mmspinorProp} The following formulas define an irreducible 
$C\!\ell_{m,m}$-module on $S = \wedge U$: 
\[ \begin{array}{r@{\: =\:}l}
\rho (u)s & u\wedge s\\
\rho (u^{\ast})s & -u^{\ast} \angle s\, , \: s\in \wedge U\, ,\: 
u\in U\, ,\: u^{\ast}\in U^{\ast}\, .
\end{array} \] 
\end{prop}

\noindent
{\bf Proof:} This follows from the obvious identities 
$\rho (u)^2 = \rho (u^{\ast})^2 = 0$ and $\rho (u) \rho (u^{\ast}) 
+ \rho (u^{\ast})\rho (u) = -2{<u,u^{\ast}>} Id$. $\Box$

For any $a\in \wedge U$ and $\alpha \in \wedge U^{\ast}$ we define 
nilpotent endomorphisms $\epsilon_a$ and $\iota_{\alpha}$ of 
$S = \wedge U$ by: 
\[ \begin{array}{r@{\: =\:}l}
\epsilon_a & a\wedge s\, ,\\
\iota_{\alpha} & \alpha \angle s\, .
\end{array}  \]

\begin{prop} The Lie algebra 
$\fr{so}(m,m)\hookrightarrow End\, S$ of the spinor group admits the
 following graded decomposition: 
\[ \fr{so}(m,m)  \, = \, \fr{g}^{-2} + \fr{g}^0 + \fr{g}^{2} \, = \, 
\iota_{\wedge^2U^{\ast}} + \fr{sl}(U) + \epsilon_{\wedge^2U}\, ,\] 
$\fr{sl}(U) = [\iota_{U^{\ast}},\epsilon_U]$, $[\fr{g}^i,\fr{g}^j] \subset 
\fr{g}^{i+j}$ ($\fr{g}^{i+j} = 0$ for $|i+j|>2$). In particular, 
$\iota_{\wedge^2U^{\ast}}$ and $\epsilon_{\wedge^2U}$ are Abelian 
subalgebras. 
\end{prop}

It is very easy to describe the semi spinor modules $S^{\pm}$ in our 
model of the spinor module $S$. 
\begin{lemma} $S = \wedge U$ is the sum of the two inequivalent irreducible
$\fr{so}(m,m)$-submodules $S^+ = \wedge^{ev}U$ and $S^- =\wedge^{odd}U$. 
\end{lemma} 

\noindent
{\bf Proof:} It is clear that $\wedge^{ev}U$ and $\wedge^{odd}U$ are 
irreducible $\fr{so}(m,m)$-submodules and we already know that they are 
inequivalent, s.\ e.g.\ Proposition \ref{redProp}. $\Box$ 

\noindent 
{\bf Remark 4:} The statement that $\wedge^{ev}U$ and $\wedge^{odd}U$ are 
inequivalent $\fr{so}(m,m)$-modules follows also from the fact that these
are  eigenspaces of the volume element $\omega_{m,m}
=e_1\cdots e_{2m}\in C\!\ell^0_{m,m}$, $(e_i)$ an orthonormal basis of 
${\Bbb R}^{m,m}$.

We can define an $\fr{so}(m,m)$-invariant endomorphism $E$ of $S$ by 
\[E|S^{\pm} = {\pm} {Id}\, . \] 

\noindent 
To construct an admissible bilinear form $f$ on $S = \wedge U$ we fix 
a volume form $vol \in \wedge^m U$ on $U^{\ast}$ and define 
\[ f(\wedge^iU,\wedge^jU)= 0\, , \quad \mbox{if} \quad i+j\neq m\, ,\] 
\[ f(s,t){vol} = \epsilon_i s\wedge t\, , \quad s\in \wedge^iU\, , \: 
t\in \wedge^{m-i}U\, ,\] 
where $\epsilon_i = (-1)^{i(i+1)/2}$. Remark that $\epsilon_{i+1} = 
(-1)^{i+1}\epsilon_i$. 
\begin{prop}\label{ffeProp} 

The space $\cal B$ of $\fr{so}(m,m)$-invariant bilinear forms on 
$S = S_{m,m}$ is 
spanned by the admissible elements $f$ and $f_E = f(E\cdot ,\cdot )$. 
Their fundamental invariants $(\tau ,\sigma , \iota )$ depend only 
on $m\pmod{4}$ and are given in the next table:\\  

\begin{centering}
\vskip6pt
\begin{tabular}{|l||c|c|c|c|}\hline  
$f$ & $- - -$ & $- - +$ & $- + -$ & $- + +$ \\ \hline 
$f_E$ & $+ + -$ & $+ - +$ & $+ - -$ & $+ + +$\\ \hline\hline
$m:$  & $1$ & $2$ & $3$ & $4$\\ \hline
\end{tabular}
\vskip6pt

\end{centering}
 An $f$- and $f_E$-admissible basis for the 
Schur algebra ${\cal C} 
\cong {\Bbb R}\oplus {\Bbb R}$ is given by the endomorphisms
$ Id$ and $E$ of $S$: 
\[ \tau (E)\: =\: -1\, ,\quad \sigma_f(E)=\sigma_{f_E}(E)\: =\: 
(-1)^m\, ,\quad 
\iota (E) = +1\, .\] 
\end{prop}

\noindent 
{\bf Proof:} We first check that $\rho (v)$, $v\in U + U^{\ast}$, is 
$f$-skew symmetric. For $v = u\in U$, $s\in \wedge^iU$, 
$t\in \wedge^{m-i-1}U$: 
\[ (f(\rho (u)s,t) + f(s,\rho (u)t)){vol} = 
\epsilon_{i+1}(u\wedge s)\wedge t + \epsilon_i s\wedge (u\wedge t) = 0\, .\] 
For $v = u^{\ast} \in U^{\ast}$, $s\in \wedge^iU$, $t\in \wedge^{m-i+1}U$: 
\begin{eqnarray*}  (f(\rho (u^{\ast})s,t) + f(s,\rho (u^{\ast})t)){vol} &=&  
\epsilon_{i-1}(u^{\ast}\angle s)\wedge t + \epsilon_is\wedge 
(u^{\ast}\angle t) = 
\end{eqnarray*} 
\[  \epsilon_{i-1}(u^{\ast}\angle s)\wedge t + 
\epsilon_i(-1)^i(u^{\ast}\angle (s\wedge t) - (u^{\ast}\angle s)\wedge t) =\]  
\[ (\epsilon_{i-1}-(-1)^i\epsilon_i)(u^{\ast}\angle s)\wedge t = 0\, .\] 
 
\noindent 
The symmetry properties of $f$ follow from the computation 
\[ f(t,s){vol} = \epsilon_jt\wedge s = \epsilon_j \epsilon_i(-1)^{ij}
f(s,t) {vol} = (-1)^{m(m+1)/2}f(s,t) {vol}\, ,\] 
where $s\in \wedge^iU$, $t\in \wedge^jU$ and $i+j=m$. 

\noindent
Finally, $f(\wedge^{ev}U,\wedge^{odd}U) = 0$ if $m$ is even and 
$f(\wedge^{ev}U,\wedge^{ev}) = f(\wedge^{odd}U,\wedge^{odd}U)$ $= 0$ 
if $m$ is odd. This proves all the statements about $f$. 
It is immediate to see that $E$ is $f$-admissible with fundamental 
invariants given above. Since $f$ is admissible and $E$ is $f$-admissible, 
$f_E$ is admissible and its fundamental invariants are computed by 
multiplicativity: 
\[ \tau (f_E) = \tau (f) \tau (E)\, , \quad \sigma (f_E) = 
\sigma (f)\sigma_f(E)\, ,\quad \iota (f_E) = \iota  (f) \iota  (E)\, .\]
This proves the proposition. $\Box$\\ 
Proposition \ref{ffeProp} implies the following theorem: 
\begin{thm} Every $\fr{so}(m,m)$-equivariant embedding 
$V^{\ast}\hookrightarrow S^{\ast}\otimes  S^{\ast}$, $S  = S_{m,m}$ the 
spinor $\fr{so}(m,m)$-module, is a linear combination of the embeddings 
$j_{\rho}(f)$ and $j_{\rho}(f_E)$. Their image is contained in the dual
of the subspaces indicated in the table depending on $m\pmod{4}$. \\   

\begin{centering}
\vskip6pt
\begin{tabular}{|l||c|c|c|c|}\hline 
$j_{\rho}(f)$ & $\vee^2S^+ + \vee^2S^-$ & $S^+\vee S^-$ & 
$\wedge^2S^+ + \wedge^2S^-$ & $S^+\wedge S^-$\\ \hline 
$j_{\rho}(f_E)$ & $\vee^2S^+ + \vee^2S^-$ & $S^+\wedge S^-$ & 
$\wedge^2S^+ + \wedge^2S^-$ & $S^+\vee S^-$ \\ \hline\hline 
$m$ & $1$ & $2$ & $3$ & $4$\\\hline  
\end{tabular}

\vskip6pt
\end{centering}
\end{thm} 
  
Now put $V_1 = {\Bbb R}^{m,m}\neq 0$ and let $V_2$ be an arbitrary 
pseudo Euclidean space. Denote the spinor module of $\fr{so}(V_i)$ 
by $S_i$, $i=1,2$. 
\begin{prop}\label{b1b2Prop} Let $\beta_2$ be 
an admissible bilinear form on $S_2$. 
Then there is a unique (up to scaling) admissible form $\beta_1$ on $S_1$ 
such that $\tau (\beta_2) = \iota (\beta_1)\tau (\beta_1)$. In particular, 
$\beta_1 \otimes \beta_2$ is an admissible bilinear form on the spinor 
$\fr{so}(V_1 + V_2)$-module $S_1 \otimes S_2$. 

If moreover, $A_2$ is a $\beta_2$-admissible endomorphism of $S_2$, then
there is a unique $\beta_1$-admissible endomorphism $A_1$ of $S_1$ such 
that $\tau (A_2) = \iota (A_1)\tau (A_1)$, in partiular, $A_1\otimes A_2$ is a 
$\beta_1 \otimes \beta_2$-admissible endomorphism of $S_1 \otimes S_2$.

The fundamental invariants of $\beta_1 \otimes \beta_2$ and 
$A_1\otimes A_2$ are easily computed using the rules given in 
Proposition \ref{multiplProp}. 
\end{prop}

\noindent
{\bf Proof:} This follows from $\iota (f_E)\tau (f_E) = 
- \iota (f) \tau (f)$, $\iota (E)\tau (E) = 
- \iota (Id) \tau (Id)$ and section \ref{redbasicsignSec}. $\Box$ 

If we assume that $V_2$ is of definite signature, i.e.\ $V_2 = 
{\Bbb R}^{k,0}$ or ${\Bbb R}^{0,k}$, then there is a 
unique (up to scaling) ${Pin} (V_2)$-invariant symmetric bilinear
form $h_2$ on the irreducible module $S_2$ of the compact group
${Pin} (V_2)$. 
\begin{lemma} The ${Pin} (V_2)$-invariant scalar product $h_2$ is admissible:
$\tau (h_2) = -1$ if $V_2 = {\Bbb R}^{k,0}$ and 
$\tau (h_2) = +1$ if $V_2 = {\Bbb R}^{0,k}$; 
$\sigma (h_2) = +1$ and if $S_2$ is reducible, $S_2 = S_2^+ + S_2^-$, 
$S_2^-=C\!\ell^1 (V_2) S_2^+$, then $\iota (h_2) = +1$. 
\end{lemma}

\noindent
{\bf Proof:} Let $\rho (v)$ denote Clifford multiplication by a unit vector 
$v\in V_2$. Then $h_2$ is  $\rho (v)$-invariant and $\rho (v)^2 = -{Id}$ 
if $V_2 = {\Bbb R}^{k,0}$ and $\rho (v)^2 = +{Id}$ if 
$V_2 = {\Bbb R}^{0,k}$. This implies $\tau (h_2) = \mp 1$. 

To see that  $\iota (h_2) = +1$ in the reducible case, consider the scalar 
product $h_2'$ on $S_2$ defined by 
\[ h_2'(S_2^+,S_2^-) = 0\, , \quad h_2'|S_2^{\pm} = h_2|S_2^{\pm}\: 
(\neq 0)\, .\] 
It is easy to check that $h_2'$ is invariant under Clifford multiplication by
unit vectors $v\in V_2$ using that $S^- = v S^+$. This implies 
$h_2' = h_2$. $\Box$ 

By Proposition \ref{b1b2Prop} for every $V_1 = {\Bbb R}^{m,m}\neq 0$ 
there is a unique admissible bilinear form $h_1$ on the spinor 
module $S_1$ of $\fr{so}(V_1)$ such that $\tau (h_2) = \iota (h_1)
\tau (h_1)$. 
\begin{dof}\label{canDef} The {\bf canonical bilinear form} on the 
spinor module 
$S = S_1 \otimes S_2$ of $\fr{so}(V_1 + V_2)$ is $h = h_1 \otimes h_2$,  
where $h_2$ is the canonical bilinear form  on the spinor module $S_2$ 
of $\fr{so}(V_2)\cong\fr{so}(k)$, i.e.\ the ${Pin}(V_2)$-invariant 
scalar product.  In line with this definition we say that an endomorphism
$A$ of $S$ (respectively $A_2$ of $S_2$) is {\bf admissible} of 
{\bf symmetry} $\sigma (A) = \pm 1$ (respectively $\sigma (A_2) = \pm 1$)
if $A$ is $h$-admissible (respectively $h_2$-admissible) and 
$\sigma_h(A) =\pm 1$ (respectively $\sigma_{h_2}(A_2) = \pm 1$).   
\end{dof}

\noindent 
{\bf Remark 5:} For $V_1 = {\Bbb R}^{m,m}$ we have two (non degenerate)
admissible bilinear forms $f$ and $f_E$ on $S_1= S_{m,m}$. 
If we want to choose 
a {\em canonical} one, which is not necessary for our purpose, we can 
consider on $S_1$ the structure of irreducible $C\!\ell_{m,m+1}$-module
defined in section \ref{cxcaseSec}.   Then only one of the forms, namely 
$f_E$, remains admissible for the $C\!\ell_{m,m+1}$-module $S_1=S_{m,m+1}$, 
it is in fact the canonical bilinear form on this module.  
Moreover,  the complex bilinear extension $f_E^{\Bbb C}$ of $f_E$, 
is the unique (up to scaling) $\fr{so}(2m+1,{\Bbb C})$-invariant complex 
bilinear form on the irreducible ${\Bbb C}\!\ell_{2m+1}$-module 
${\Bbb S}_{2m+1} = S_{m,m+1}\otimes{\Bbb C}$, s.\ Corollary 
\ref{oddcxCor}.

\subsection{Complex case}\label{cxcaseSec} 
{\bf Case of even dimension:} \\ 
The following theorem follows immediately from the fact that an irreducible 
module ${\Bbb S}_{2m}$ of ${\Bbb C}\!\ell_{2m}$ can be 
obtained as ${\Bbb S}_{2m} = S_{m,m}\otimes {\Bbb C}$ and that  
${\Bbb S}_{2m}$ splits as 
${\Bbb C}\!\ell^0_{2m}$-module: 
${\Bbb S}_{2m} = {\Bbb S}_{2m}^+ + {\Bbb S}_{2m}^-$, where 
${\Bbb S}_{2m}^{\pm} = S^{\pm}_{m,m}\otimes {\Bbb C}$. 

\begin{thm} Every $\fr{so}(2m,{\Bbb C})$-equivariant embedding 
${\Bbb C}^{2m} \hookrightarrow {\Bbb S}_{2m}\otimes 
{\Bbb S}_{2m}$ is a linear combination of the embeddings 
$j_{\rho}(f)^{\Bbb C}$ and $j_{\rho}(f_E)^{\Bbb C}$. 
Their image is contained in the dual
of the subspaces indicated in the table depending on $m\pmod{4}$, 
where we have put ${\Bbb S} = {\Bbb S}_{2m}$.\\  

\begin{centering}
\vskip6pt
\begin{tabular}{|l||c|c|c|c|}\hline 
$j_{\rho}(f)^{\Bbb C}$ & $\vee^2{\Bbb S}^+ + \vee^2{\Bbb S}^-$ 
& ${\Bbb S}^+\vee {\Bbb S}^-$ & 
$\wedge^2{\Bbb S}^+ + \wedge^2{\Bbb S}^-$ & ${\Bbb S}^+\wedge 
{\Bbb S}^-$\\ \hline 
$j_{\rho}(f_E)^{\Bbb C}$ & $\vee^2{\Bbb S}^+ + 
\vee^2{\Bbb S}^-$ & ${\Bbb S}^+\wedge {\Bbb S}^-$ & 
$\wedge^2{\Bbb S}^+ + \wedge^2{\Bbb S}^-$ & ${\Bbb S}^+ 
\vee {\Bbb S}^-$ \\ \hline\hline 
$m$ & $1$ & $2$ & $3$ & $4$\\\hline  
\end{tabular}

\vskip6pt
\end{centering}
\end{thm} 
 
\noindent
{\bf Case of odd dimension:} \\ 
The odd dimensional complex case can be obtained from the real case of 
signature $(m,m+1)$ by complexification. 

We fix the orthogonal decomposition $({\Bbb R}^{m,m+1},{<\cdot ,\cdot >})
= {\Bbb R}e_0 + {\Bbb R}^{m,m}$, where ${<e_0,e_0>} = -1$, and 
denote by $\rho$ the irreducible representation of $C\!\ell_{m,m}$ on 
$S_{m,m}$ constructed in Proposition \ref{mmspinorProp}. 
\begin{prop}An irreducible representation $\tilde{\rho}$ of 
$C\!\ell_{m,m+1}$ on $S_{m,m+1} = S_{m,m}$ is defined by 
\[ \tilde{\rho}|{\Bbb R}^{m,m} = \rho |{\Bbb R}^{m,m}\, ,\quad 
\tilde{\rho}(e_0) = \rho (\omega_{m,m} )\, ,\] 
where $\omega_{m,m}$ is the volume element of $C\!\ell_{m,m}$. 
The  $C\!\ell_{m,m+1}^0$-module $S_{m,m+1}$ is irreducible and has 
Schur algebra ${\cal C}_{m,m+1} = {\Bbb R}\, {Id}$.  
\end{prop} 

\noindent 
{\bf Proof:} It is sufficient to check that $\{ \tilde{\rho}(e_0), 
\rho (x)\} = 0$ for $x\in {\Bbb R}^{m,m}$ and that 
$\tilde{\rho}(e_0)^2 = {Id}$. This follows from the next lemma. $\Box$ 
\begin{lemma} The volume element $\omega = \omega_{m,m} = e_1e_2 \cdots 
e_{2m}$ ($(e_i)$ an  orthonormal basis of  
${\Bbb R}^{m,m}$) of $C\!\ell_{m,m}$ satisfies $\{ \omega ,x\} = 0$ 
for all $x\in {\Bbb R}^{m,m}$ and $\omega^2 = +1$. 
\end{lemma}

\begin{prop} Every $\fr{so}(m,m+1)$-invariant bilinear form on 
$S = S_{m,m+1}$ is a multiple of the admissible (canonical) form
$f_E$ (s.\ Proposition \ref{ffeProp}) and hence every  
$\fr{so}(m,m+1)$-equivariant embedding ${\Bbb R}^{m,m+1} 
\hookrightarrow (S\otimes S)^{\ast}$ is proportional to the 
embedding $j_{\tilde{\rho}}(f_E)$, which maps 
${\Bbb R}^{m,m+1}$ into $\vee^2 S^{\ast}$ if $m\equiv 0$ or $1\pmod{4}$ 
and into $\wedge^2S^{\ast}$ if $m\equiv 2$ or $3\pmod{4}$. 
\end{prop} 

\noindent 
{\bf Proof:} $\tilde{\rho}(e_0) = \rho (\omega_{m,m})$ is $f_E$-symmetric and 
$\tau (f_E)= +1$. $\Box$ 
    
\begin{cor}\label{oddcxCor} Every $\fr{so}(2m+1,{\Bbb C})$-invariant 
bilinear form on ${\Bbb S} = {\Bbb S}_{2m+1} = 
S_{m,m+1}\otimes {\Bbb C}$ is a multiple of the form 
$f_E^{\Bbb C}$ and every $\fr{so}(2m+1,{\Bbb C})$-equivariant 
embedding ${\Bbb C}^{2m+1} \hookrightarrow 
({\Bbb S}\otimes {\Bbb S})^{\ast}$ is proportional to the 
embedding $j_{\tilde{\rho}}(f_E)^{\Bbb C}$. 
\end{cor}
\section{Case of positive signature}\label{posSec} 
\subsection{Case of even dimension}\label{posevenSubsec}  
We fix the orthogonal decomposition ${\Bbb R}^{2m} = {\Bbb R}^m + 
\widetilde{{\Bbb R}^m}$, where 
$\, \tilde{{ }}: {\Bbb R}^m \rightarrow
\widetilde{{\Bbb R}^m}$ is an isometry. Denote by $\alpha$ the involution
of $C\!\ell_m$ (respectively ${\Bbb C}\!\ell_m$) extending $x\mapsto -x$ 
on ${\Bbb R}^m$ (respectively ${\Bbb C}^m$). 
\begin{prop} If $m \equiv 0$ or $3\pmod{4}$ the following formulas define on 
$S = S_{2m,0} = C\!\ell_m$ the structure of irreducible  
$C\!\ell_{2m}$-module: 
\begin{eqnarray*}
\rho (x)s &=& xs\\
\rho (\tilde{x})s &=& \omega sx \quad \mbox{if} \quad m\equiv 0\pmod{4}\\
\rho (\tilde{x})s &=& \omega \alpha (s) x \quad \mbox{if} \quad m\equiv 
3\pmod{4}\, ,
\end{eqnarray*} 
where $x \in {\Bbb R}^m$, $s\in S$ and $\omega$ is the volume element 
of $C\!\ell_m$, i.e.\ $\omega = e_1\cdots e_m$ for an orthonormal 
basis $(e_i)$ of  ${\Bbb R}^m$. The $\fr{so}(2m)$-module $S$ is the 
sum $S = S^+ + S^-$ of the two inequivalent irreducible modules 
$S^+ = C\!\ell_m^0$ and $S^- = C\!\ell_m^1$ if $m\equiv 0\pmod{4}$ and 
is irreducible if $m\equiv 3\pmod{4}$. 

If $m \equiv 1$ or $2\pmod{4}$ the structure of irreducible  
$C\!\ell_{2m}$-module on $S = S_{2m,0} = {\Bbb S}_{2m} = 
{\Bbb C}\!\ell_m$ is given by: 
\begin{eqnarray*}
\rho (x)s &=& xs\\
\rho (\tilde{x})s &=& i\alpha (s) x\, , \quad x\in {\Bbb R}^m\, , \quad 
s\in S\, .
\end{eqnarray*} 
As $\fr{so}(2m)$-module $S = S^+ + S^-$ is the 
sum   of the two irreducible modules 
$S^+ = {\Bbb C}\!\ell_m^0$ and 
$S^- = {\Bbb C}\!\ell_m^1$, which are equivalent for 
$m \equiv 1\pmod{4}$ and inequivalent for $m \equiv 2\pmod{4}$. 
\end{prop} 

\noindent
{\bf Proof:} It is sufficient to check the identities 
{\arraycolsep=0pt
\begin{eqnarray*}
\rho (x)^2 &{}={}& -{<x,x>} {Id},\\
\rho (\tilde{x})^2 &{}={}& -{<x,x>} {Id}, \\
\{ \rho (x),\rho (\tilde{y})\} &{}={}& 0 
\end{eqnarray*}}
for 
$x,y \in {\Bbb R}^m$. This is straightforward using the following 
lemma. $\Box$ 

\begin{lemma} \label{wmLemma} The volume element $\omega = \omega_m = e_1 \cdots e_m$ 
of $C\!\ell_m$ satisfies $\{ \omega,x\} = 0$ if m is even and 
$[\omega ,x] = 0$ if m is odd, $x\in {\Bbb R}^m \subset C\!\ell_m$. 
Moreover, 
\[ \omega^2 \: = \: \left\{ \begin{array}{r@{\quad\mbox{if}\quad m\equiv\:}l} 
+1 & 0\quad \mbox{or} \quad 3\pmod{4}\\
-1 & 1\quad \mbox{or} \quad 2\pmod{4}\, .
\end{array}\right.  \] 
\end{lemma} 
 
Now we describe the ${Pin} (2m)$-invariant symmetric bilinear form $h$ on $S$ 
using the canonical identification $\wedge {\Bbb R}^m \rightarrow 
C\!\ell_m$ of {\mbox Z}$_2$-graded vector spaces given by 
\[ e_{i_1}\wedge \ldots \wedge e_{i_k} \mapsto e_{i_1}\cdots e_{i_k}\] 
with respect to an orthonormal basis $(e_i)$, $i=1,\ldots ,m$, of 
${\Bbb R}^m$. 

The standard scalar product $<\cdot ,\cdot >$ on $\wedge {\Bbb R}^m$ 
induced by the scalar product on ${\Bbb R}^m$ is invariant under 
exterior  $x\wedge \cdot$  and interior $x\angle \cdot$ multiplication 
with unit vectors 
$x \in {\Bbb R}^m$. 
\begin{lemma} Using the identification $C\!\ell_m = \wedge {\Bbb R}^m$, 
Clifford multiplication of $x\in {\Bbb R}^m$ and $\phi \in C\!\ell_m$
is given by:
\begin{eqnarray*} 
x\phi &=& x\wedge \phi - x\angle  \phi\\ 
\phi x &=& x\wedge \alpha (\phi ) + x\angle  \alpha (\phi )\, .
\end{eqnarray*} 
\end{lemma} 

\noindent 
{\bf Proof:} The proof is similar to \cite{L-M} I.Prop.\ 3.9. $\Box$ 

\begin{cor} \label{invarianceCor} The standard scalar 
product $<\cdot , \cdot >$ on 
$\wedge {\Bbb R}^m = C\!\ell_m$ is invariant under left  and right
multiplications by unit vectors   $x \in {\Bbb R}^m$. In particular, 
if $m\equiv 0$ or $3\pmod{4}$, $h = {<\cdot ,\cdot >}$ is the 
(admissible) ${Pin} (2m)$-invariant  scalar product on the irreducible 
$C\!\ell_{2m}$-module $S = C\!\ell_m$. 
\end{cor} 

If $m \equiv 1$ or $2\pmod{4}$, we extend the standard scalar product  on 
$\wedge {\Bbb R}^m$ to a symmetric complex bilinear form 
${<\cdot ,\cdot >}_{{\Bbb C}}$ on $S = \wedge {\Bbb C}^m$. 
Using the operator $c$ of complex conjugation, we define a symmetric 
real bilinear form $h = Re\, {<c\, \cdot ,\cdot >}_{{\Bbb C}}$ on $S$. 

\begin{lemma} Let $m \equiv 1$ or $2\pmod{4}$. Then  $h = 
Re\, {<c\, \cdot ,\cdot >}_{{\Bbb C}}$ is the (admissible) 
${Pin} (2m)$-invariant  scalar product on the irreducible 
$C\!\ell_{2m}$-module $S = {\Bbb C}\!\ell_m$. 
\end{lemma} 

\noindent
{\bf Proof:}  We check that $\rho (x)$ and $\rho (\tilde{x})$, $x\in 
{\Bbb R}^m$, are ${<c\, \cdot ,\cdot >}_{{\Bbb C}}$-skew symmetric 
and hence $h$-skew symmetric. By Corollary \ref{invarianceCor} left 
and right multiplication, $L_x$ and $R_x$, by $x\in 
{\Bbb R}^m$ are ${<\cdot ,\cdot >}_{{\Bbb C}}$-skew symmetric
endomorphisms of $S = {\Bbb C}\!\ell_m$, in particular,  
$\rho (x)$ is ${<\cdot ,\cdot >}_{{\Bbb C}}$-skew symmetric. 
It is easy to see that $\alpha$ and the operator $I$ of multiplication 
by $i$ are ${<\cdot ,\cdot >}_{{\Bbb C}}$-symmetric endomorphisms. 
Moreover, 
\[ [I,R_x] = [I,\alpha ] = \{ \alpha ,R_x\} = 0\] 
and hence $\rho (\tilde{x}) = I\circ R_x \circ \alpha$ is 
${<\cdot ,\cdot >}_{{\Bbb C}}$-symmetric. From the relations 
\[ [c,L_x] = [c,R_x] = [c,\alpha ] = \{ c,I\} = 0\] 
we obtain that $[\rho (x),c] = \{ \rho (\tilde{x}),c\} = 0$, 
which implies that  $\rho (x)$ and $\rho (\tilde{x})$ are 
${<c\, \cdot ,\cdot >}_{{\Bbb C}}$-skew symmetric. $\Box$ 

Now we construct admissible, i.e.\ $h$-admissible, bases  of  
the Schur algebra ${\cal C} = {\cal C}_{2m,0}$ for all the values of 
$m\pmod{4}$. 

\begin{prop} \label{2mSchurProp} If $m \equiv 0\pmod{4}$,  an admissible 
basis of the Schur 
algebra 
$${\cal C}_{2m,0} \cong {\Bbb R} \oplus {\Bbb R}$$ is 
given by the endomorphisms $ Id$ and $E=\alpha$ of $S = C\!\ell_m$: 
$\tau (E) = -1$, $\sigma (E) = \sigma_h(E) = +1$, $\iota (E) = +1$. 

If $m \equiv 3\pmod{4}$,  an admissible basis of 
${\cal C}_{2m,0} \cong {\Bbb C}$ is given by the endomorphisms 
$ Id $ and
$J = L_{\omega} \circ \alpha$ of $S = C\!\ell_m$: $\tau (J) = -1$, $\sigma (J) 
= -1$. 

The space $\cal B$ of $\fr{so}(2m)$-invariant bilinear forms on $S$ is 
spanned by admissible elements: 
\[ {\cal B} = span \, \{ h,h_E\} \quad \mbox{if}\quad m \equiv 0\pmod{4}\, ,\] 
\[ {\cal B} = span \, \{ h,h_J\} \quad \mbox{if}\quad m \equiv 3\pmod{4}\, .\]
The fundamental invariants $(\tau ,\sigma , \iota )$ are given by:
\begin{eqnarray*}
(\tau ,\sigma , \iota )(h) &=& (-1,+1,+1), \\
(\tau ,\sigma , \iota )(h_E)&=&(+1,+1,+1)\mbox{ if }m \equiv 0\pmod{4},\\
(\tau ,\sigma )(h) &=& (-1,+1), \\
(\tau ,\sigma )(h_J) &=& (+1,-1)\mbox{ if }m \equiv 3\pmod{4}. 
\end{eqnarray*}
\end{prop}

\noindent 
{\bf Proof:} We show that $J$ is admissible and $\tau (J) = \sigma (J) = -1$. 
All other statements are immediate. 

Let $m \equiv 3\pmod{4}$. 
From $[L_x,L_{\omega}] = [R_x,L_{\omega}] = \{ L_x, \alpha \} = \{ R_x, 
\alpha \} = 0$   (s.\ Lemma \ref{wmLemma}) it follows that 
$\{ L_x,J\} = \{ R_x,J\} = 0$. Since $\rho (x) = L_x$ and 
$\rho (\tilde{x}) = R_x\circ J$, we conclude 
$\{ \rho (x),J\} = \{ \rho (\tilde{x}), J\} = 0$. \\ 
The operator $J$ is skew symmetric as  product of  two anticommuting 
symmetric  operators, namely $L_{\omega}$ and $\alpha$ (the scalar 
product is $L_{\omega}$-invariant and $L_{\omega}^2 = + {Id}$). $\Box$ 

If $m \equiv 1$ or $2\pmod{4}$, we consider the following operators on $S
= {\Bbb C}\!\ell_m$: 
\[ I: s\mapsto is\, ,\: J = L_{\omega}\circ c\, ,\: K= IJ \quad \mbox{and} 
\quad E =  \alpha \, ,\] 
where $\omega = e_1\cdots e_m\in C\!\ell_m\subset {\Bbb C}\!\ell_m$ 
is the volume element. 
\begin{prop} Let $m \equiv 1$ or $2\pmod{4}$. The Schur algebra 
${\cal C}_{2m,0}$ ($\cong {\Bbb C}(2)$ if $m \equiv 1\pmod{4}$ and 
$\cong {\Bbb H}\oplus {\Bbb H}$ if $m \equiv 2\pmod{4}$) 
is generated by the admissible operators $I$, $J$ and $E$ satisfying the 
following (anti) commutator relations: 
\[ I^2 = J^2 = L_{\omega}^2 = -1\, ,\quad E^2 =c^2 = +1\, ,\] 
\[ \{ I,J\} = [I,E] = [I,L_{\omega}] = \{I,c\} = 0\, ,\] 
\[ [J, L_{\omega}] = [J,c] =[E,c] = [L_{\omega},c] = 0\, ,\] 
\[ \{ J,E\} = \{ L_{\omega},E\} = 0 \quad \mbox{if} \quad m \equiv 1\pmod{4}
\, ,\] 
\[ [J,E] = [L_{\omega},E] = 0 \quad \mbox{if} \quad m \equiv 2\pmod{4}\, .\] 
An admissible basis of the Schur algebra is given by the endomorphisms 
$ Id $, $I$, $J$, $K$, $E$, $EI$, $EJ$, $EK$.  Their fundamental 
invariants $(\tau ,\sigma , \iota )$ are given in the next table, where 
the value of $m$ is modulo 4. \\ 

\begin{centering}
\begin{tabular}{|l||c|c|c|c|c|c|c|c|}\hline 
$m$: & $ Id $ & $I$ & $J$ & $K$ & $E$ & $EI$ & $EJ$ & $EK$\\\hline\hline
$1$ & $+++$ & $+-+$ & $+--$ & $+--$ & $-++$ & $--+$ & $-+-$ & $-+-$\\\hline 
$2$ & $+++$ & $+-+$ & $--+$ & $--+$ & $-++$ & $--+$ & $+-+$ & $+-+$\\\hline 
\end{tabular}
\vskip 6pt

\end{centering}

The fundamental invariants  of the corresponding admissible basis of $\cal B$ 
are also listed for convenience:

\begin{centering}
\hskip 6pt
\begin{tabular}{|l||c|c|c|c|c|c|c|c|}\hline 
$m$: & $h$ & $h_I$ & $h_J$ & $h_K$ & $h_E$ & $h_{EI}$ & $h_{EJ}$ & 
$h_{EK}$\\\hline\hline
$1$ & $-++$ & $--+$ & $---$ & $---$ & $+++$ & $+-+$ & $++-$ & $++-$\\\hline 
$2$ & $-++$ & $--+$ & $+-+$ & $+-+$ & $+++$ & $+-+$ & $--+$ & $--+$\\\hline 
\end{tabular} 

\end{centering}
\end{prop} 

\noindent
{\bf Proof:} The proof is similar to the proof of Proposition \ref{ffeProp} 
and \ref{2mSchurProp}. One uses the multiplication rules for the invariants 
and also that $L_{\omega}$ is skew symmetric, $c$ is symmetric and 
they commute. $\Box$ 

\begin{thm} Every $\fr{so}(2m)$-equivariant  embedding 
${\Bbb R}^{2m} \hookrightarrow (S\otimes S)^{\ast}$,
$S = S_{2m,0}$, is a 
linear combination of the embeddings 
\[ j_{\rho}(h): {\Bbb R}^{2m} \hookrightarrow (S^+\wedge S^-)^{\ast} 
\quad \mbox{and} \quad j_{\rho}(h_E): {\Bbb R}^{2m} \hookrightarrow (S^+
\vee S^-)^{\ast}\] 
if $m\equiv 0\pmod{4}$ and a 
linear combination of 
\[ j_{\rho}(h): {\Bbb R}^{2m} \hookrightarrow \wedge^2S^{\ast} 
\quad \mbox{and} \quad j_{\rho}(h_J): {\Bbb R}^{2m} \hookrightarrow 
\wedge^2S^{\ast} \] 
if $m\equiv 3\pmod{4}$.  

If $m\equiv 1$ or $2\pmod{4}$ every $\fr{so}(2m)$-equivariant  embedding 
${\Bbb R}^{2m} \hookrightarrow (S\otimes S)^{\ast}$ is a 
linear combination of the embeddings $j_A = j_{\rho}(h_A)$, $A\in {\cal C}$ 
admissible, whose image is contained in the dual of the subspaces indicated 
in Table \ref{2mTable} depending on $m\pmod{4}$. 
\begin{table}[ht]\caption[$\fr{so}(2m)$-equivariant  embeddings 
$j_A = j_{\rho}(h_A): {\Bbb R}^{2m} 
\hookrightarrow (S\otimes S)^{\ast}$]{\label{2mTable} 
$\fr{so}(2m)$-equivariant  embeddings 
$j_A = j_{\rho}(h_A): {\Bbb R}^{2m} 
\hookrightarrow (S\otimes S)^{\ast}$}  

\begin{centering}
\begin{tabular}{|l||c|c|} \hline 
$j_{Id}$ & $S^+\wedge S^-$ & $S^+\wedge S^-$\\\hline 
$j_I$ & $S^+\vee S^-$ & $S^+\vee S^-$\\\hline 
$j_J$ & $\vee^2S^+ + \vee^2S^-$ & $S^+\wedge S^-$\\\hline 
$j_K$ & $\vee^2S^+ + \vee^2S^-$ & $S^+\wedge S^-$\\\hline 
$j_E$ & $S^+\vee S^-$ & $S^+\vee S^-$\\\hline 
$j_{EI}$ & $S^+\wedge S^-$ & $S^+\wedge S^-$\\\hline 
$j_{EJ}$ & $\vee^2S^+ + \vee^2S^-$ & $S^+\vee S^-$\\\hline 
$j_{EK}$ & $\vee^2S^+ + \vee^2S^-$ & $S^+\vee S^-$\\\hline\hline 
$m$: & $1$ & $2$ \\\hline 
\end{tabular} 

\end{centering}
\vskip18pt
\end{table}  
\end{thm} 
\subsection{Case of odd dimension} 
To reduce the odd dimensional case to the even dimensional, we consider the 
orthogonal decomposition ${\Bbb R}^{2m+1} = {\Bbb R}e_0 + 
{\Bbb R}^{2m}$, where $e_0$ is a unit vector. Let $\rho$ 
denote the irreducible representation of $C\!\ell_{2m}$ on 
$S_{2m,0}$ defined in section \ref{posevenSubsec}. We  will extend
$\rho$ to an irreducible representation $\tilde{\rho}$ of $C\!\ell_{2m+1}$ 
on $S = S_{2m+1,0}$, where $S_{2m+1,0} = S_{2m,0}$ if $m\equiv 1$, $2$ 
or $3\pmod{4}$ and  $S_{2m+1,0} = S_{2m,0}\otimes {\Bbb C} = 
{\Bbb S}_{2m}$ if $m\equiv 0\pmod{4}$. If $m\equiv 1$ or 
$2\pmod{4}$, $S_{2m,0} = {\Bbb S}_{2m}$ admits the 
$C\!\ell_{2m}$-invariant complex structure $I$. For $m\equiv 0\pmod{4}$ 
multiplication by $i$ is a $C\!\ell_{2m}$-invariant complex structure
on $S_{2m,0}\otimes {\Bbb C}$ and will also be denoted by $I$. 

\begin{prop} The following formulas define an irreducible representation $\tilde{\rho}$ of $C\!\ell_{2m+1}$ 
on $S_{2m+1,0}$. 
\[ \tilde{\rho}|{\Bbb R}^{2m} = \rho |{\Bbb R}^{2m}\, ,\] 
\[ \tilde{\rho}(e_0) = \left\{ 
\begin{array}{r@{\quad \mbox{if} \quad m\equiv \,}l} 
\rho (\omega_{2m}) & 1\quad \mbox{or} \quad 3\pmod{4} \\
I \circ \rho (\omega_{2m}) & 0\quad \mbox{or} \quad 2\pmod{4}\, ,
\end{array} \right. \] 
where, in the case   $m\equiv 0\pmod{4}$,  $\rho$ has been extended 
complex linearly to a representation  on $S_{2m,0}\otimes {\Bbb C}$, 
denoted by the same symbol. $S =  S_{2m+1,0}$ is irreducible as 
$C\!\ell_{2m+1}^0$-module if $m\not\equiv 0\pmod{4}$ and the sum 
$S = S^+ + S^-$ of the two equivalent irreducible  $C\!\ell_{2m+1}^0$-modules 
$S^+ = S^+_{2m,0} + iS^-_{2m,0} = C\!\ell_m^0 + iC\!\ell_m^1$ and 
$S^- = iS^+$ if $m\equiv 0\pmod{4}$. 
\end{prop} 

\noindent 
{\bf Proof:} It is sufficient to check that $\tilde{\rho}(e_0)^2 = -{Id}$ and 
$\{ \tilde{\rho}(e_0), \rho (x) \} = 0$ for  $x\in {\Bbb R}^{2m}$, 
since all other information can be extracted from the Schur algebra, 
s.\ Corollary \ref{SchurCor}. These identities follow immediately from 
Lemma \ref{wmLemma} and the fact that $I$ is a $C\!\ell_{2m}$-invariant 
complex structure. $\Box$ 

Now we describe the ${Pin} (2m+1)$-invariant scalar product $h$ on 
$S =  S_{2m+1,0}$. Let $h_{2m,0}$ denote the ${Pin} (2m)$-invariant 
scalar product  on $S_{2m+1,0} = S_{2m,0}$ if $m\equiv 1$, $2$ or 
$3\pmod{4}$ and by $h_{2m,0}^{\Bbb C}$ the complex bilinear 
extension of the  ${Pin} (2m)$-invariant 
scalar product  on $S_{2m,0}$  to a ${Pin} (2m)$-invariant complex 
bilinear form on $S_{2m+1,0} = {\Bbb S}_{2m} =  S_{2m,0}
\otimes {\Bbb C}$ if $m\equiv 4\pmod{4}$. 
\begin{lemma} The ${Pin} (2m+1)$-invariant scalar product $h = h_{2m+1,0}$ 
on $S =  S_{2m+1,0}$ is given by $h = h_{2m,0}$ if $m\equiv 1$, $2$ or 
$3\pmod{4}$ and by $h = Re\, h_{2m,0}^{\Bbb C}(c\, \cdot ,\cdot )$ if 
$m\equiv 4\pmod{4}$, where $c$ is complex conjugation with respect to 
$S_{2m,0}\subset S_{2m,0}\otimes {\Bbb C}$. 
\end{lemma} 

\noindent 
{\bf Proof:} If $m\not\equiv 4\pmod{4}$, the statement follows from Schur's 
Lemma, since  $S_{2m+1,0} = S_{2m,0}$. If $m\equiv 4\pmod{4}$, 
the Hermitian form $h_{2m,0}^{\Bbb C}(c\, \cdot ,\cdot )$ is 
$I$-invariant and hence invariant under $\tilde{\rho}(e_0) = 
I \circ \rho (\omega_{2m})$ and the same is true for 
$h = Re\, h_{2m,0}^{\Bbb C}(c\, \cdot ,\cdot )$. $\Box$ 

If $m\not\equiv 3\pmod{4}$, we have on $S_{2m+1,0} = {\Bbb C}\!\ell_m = 
C\!\ell_m + iC\!\ell_m$ the operator $c$ of complex conjugation. Hence, 
we can define an endomorphism $J$ of $S_{2m+1,0} = {\Bbb C}\!\ell_m$ 
by the formulas 
\[ J: = \left\{ \begin{array}{r@{\quad \mbox{if} \quad m\equiv \:}l}  
L_{\omega}\circ c & 1\quad \mbox{or} \quad 2\pmod{4}\\ 
\alpha \circ c & 0\pmod{4}\, ,
\end{array} \right. \]    
where $L_{\omega}$ is left multiplication by the volume element 
$\omega = \omega_m$ of $C\!\ell_m$ and $\alpha | 
{\Bbb C}\!\ell_m^0 = +{Id}$, $\alpha | 
{\Bbb C}\!\ell_m^1 = -{Id}$. 

\begin{prop} Let $m\not\equiv 3\pmod{4}$. An admissible basis of the Schur 
algebra ${\cal C} = {\cal C}_{2m+1,0}$ is given by the endomorphisms 
$ Id $, $I$, $J$ and  $K=IJ$ of $S_{2m+1,0}= {\Bbb C}\!\ell_m$. 
If $m\equiv 1$ or $2\pmod{4}$, then $I^2 = J^2 = -{Id}$, $\{ I,J\} = 0$ 
and ${\cal C}_{2m+1,0}\cong {\Bbb H}$. If $m\equiv 0\pmod{4}$, then
$I^2 = -J^2 = -{Id}$, $\{ I,J\} = 0$ 
and ${\cal C}_{2m+1,0}\cong {\Bbb R}(2)$. The space $\cal B$ of 
$\fr{so}(2m+1)$-invariant bilinear forms on $S_{2m+1,0}$ has the admissible 
basis $(h,h_I,h_J,h_K)$. If $m\equiv 3\pmod{4}$, then the Schur algebra 
${\cal C}_{2m+1,0} = {\Bbb R}\, {Id}$ and ${\cal B} = {\Bbb R}h$. 
\end{prop} 

\noindent
{\bf Proof:} straightforward, cf.\ Proposition \ref{2mSchurProp}. $\Box$ 

\begin{thm} If $m\equiv 3\pmod{4}$, every $\fr{so}(2m+1)$-equivariant  
embedding $${\Bbb R}^{2m+1} \hookrightarrow S^{\ast}\otimes S^{\ast},$$ 
$S = S_{2m+1,0}$, is a multiple of $j_{\rho}(h): {\Bbb R}^{2m+1} 
\hookrightarrow \wedge^2S^{\ast}$. If $m\not\equiv 3\pmod{4}$, every 
$\fr{so}(2m+1)$-equivariant  
embedding $${\Bbb R}^{2m+1} \hookrightarrow (S\otimes S)^{\ast}$$
is a linear combination of the embeddings $j_A = j_{\rho}(h_A)$, $A = 
{Id}$, $I$, $J$ or $K$, whose  image is contained in the dual 
of the subspaces indicated in Table \ref{2m+1Table} depending on $m\pmod{4}$.
\begin{table}[ht]\caption[$\fr{so}(2m+1)$-equivariant  embeddings 
$j_A = j_{\rho}(h_A): {\Bbb R}^{2m+1} 
\hookrightarrow (S\otimes S)^{\ast}$]{\label{2m+1Table} 
$\fr{so}(2m+1)$-equivariant  embeddings 
$j_A: {\Bbb R}^{2m+1} 
\hookrightarrow (S\otimes S)^{\ast}$} 

\begin{centering}
\begin{tabular}{|l||c|c|c|c|} \hline
$m$: & $j_{Id}$ & $j_I$ & $j_J$ &  $j_K$ \\\hline\hline 
$1$ & $\wedge^2S$ & $\vee^2S$ & $\vee^2S$ & $\vee^2S$\\\hline 
$2$ & $\wedge^2S$ & $\vee^2S$ & $\wedge^2S$ & $\wedge^2S$\\\hline 
$4$ & $S^+\wedge S^-$ & $\vee^2S^+ + \vee^2S^-$ & $S^+\vee S^-$ & 
$\vee^2S^+ + \vee^2S^-$\\\hline 
\end{tabular} 

\end{centering}
\vskip18pt
\end{table}  
\end{thm} 
\section{Case of negative signature}\label{negSec} Now we discuss the case of 
negative signature. 
The proofs are similar to the proofs in the case of positive signature and 
will mostly be omitted. 
\subsection{Case of even dimension}\label{negevenSubsec}  
As in the positively defined case, we fix the orthogonal 
decomposition ${\Bbb R}^{0,2m} = {\Bbb R}^{0,m} 
+ \widetilde{{\Bbb R}^{0,m}}$, where $\tilde{{ }}: {\Bbb R}^{0,m} 
\rightarrow \widetilde{{\Bbb R}^{0,m}}$ is an isometry. 
\begin{lemma} \label{w0mLemma} The volume element $\omega = \omega_{0,m} 
= e_1 \cdots e_m$  ($(e_i)$ an orthonormal basis of ${\Bbb R}^{0,m}$) 
of $C\!\ell_{0,m}$ satisfies $\{ \omega,x\} = 0$ if m is even and 
$[\omega ,x] = 0$ if m is odd, 
$x\in {\Bbb R}^{0,m}\subset C\!\ell_{0,m}$. Moreover, 
\[ \omega^2 \: = \: \left\{ \begin{array}{r@{\quad\mbox{if}\quad m\equiv \:}l} 
+1 & 0\quad \mbox{or} \quad 1\pmod{4}\\
-1 & 2\quad \mbox{or} \quad 3\pmod{4}\, .
\end{array}\right.  \] 
\end{lemma}  
The next proposition is checked using Lemma \ref{w0mLemma}. 
\begin{prop} If $m \equiv 0$ or $1\pmod{4}$ the following formulas define on 
$S = S_{0,2m} = C\!\ell_{0,m}$ the structure of irreducible  
$C\!\ell_{0,2m}$-module: 
\begin{eqnarray*}
\rho (x)s &=& xs\\
\rho (\tilde{x})s &=& \omega sx \quad \mbox{if} \quad m\equiv 0\pmod{4}\\
\rho (\tilde{x})s &=& \omega \alpha (s) x \quad \mbox{if} \quad m\equiv 1\pmod{4}\, ,
\end{eqnarray*} 
where $x \in {\Bbb R}^{0,m}$, $s\in S$ and $\omega$ is the volume element 
of $C\!\ell_{0,m}$. The $\fr{so}(0,2m)$-module $S$ is the 
sum $S = S^+ + S^-$ of the two inequivalent irreducible modules 
$S^+ = C\!\ell_{0,m}^0$ and $S^- = C\!\ell_{0,m}^1$ if $m\equiv 0\pmod{4}$ and 
is irreducible if $m\equiv 1\pmod{4}$. 

If $m \equiv 2$ or $3\pmod{4}$ the structure of irreducible  
$C\!\ell_{0,2m}$-module on $S = S_{0,2m} = {\Bbb S}_{2m} = 
{\Bbb C}\!\ell_m$ is given by: 
\begin{eqnarray*}
\rho (x)s &=& xs\\
\rho (\tilde{x})s &=& i\alpha (s) x\, , \quad x\in {\Bbb R}^{0,m}
\subset {\Bbb C}\!\ell_m = C\!\ell_{0,m}\otimes {\Bbb C}\, , 
\: s\in S = {\Bbb C}\!\ell_m\, .
\end{eqnarray*} 
As $\fr{so}(0,2m)$-module $S = S^+ + S^-$ is the 
sum   of the two irreducible sub\-mo\-dules 
$S^+ = {\Bbb C}\!\ell_m^0$ and 
$S^- = {\Bbb C}\!\ell_m^1$, which are inequivalent for 
$m \equiv 2\pmod{4}$ and equivalent for $m \equiv 3\pmod{4}$. 
\end{prop} 

Recall (s.\ Corollary \ref{invarianceCor}) that the standard scalar 
product  on 
$\wedge {\Bbb R}^m = C\!\ell_m =  C\!\ell_{m,0}$ is invariant under
left and right multiplications by unit vectors   $x \in {\Bbb R}^m = 
{\Bbb R}^{m,0}$. We can consider ${\Bbb R}^{0,m}$ as subspace 
\[ {\Bbb R}^{0,m} = i{\Bbb R}^m\subset {\Bbb C}\!\ell_m =
C\!\ell_m \otimes  {\Bbb C} = C\!\ell_m + iC\!\ell_m\, .\] 
Then $C\!\ell_{0,m} = C\!\ell_{0,m}^0 + C\!\ell_{0,m}^1 = 
C\!\ell_m^0 + iC\!\ell_m^1$. We define an isomorphism of 
$\Bbb Z_2$-graded vector spaces $\varphi : C\!\ell_m\rightarrow 
C\!\ell_{0,m}$ on elements $a\in C\!\ell_m$ of pure degree $\deg (a) = 0$ or 
$1$ by: 
\[ a\mapsto i^{\deg (a)} a\, .\] 
A  scalar product $<\cdot , \cdot >$ on  $C\!\ell_{0,m}$ is defined by the 
condition that $\varphi : C\!\ell_m\rightarrow 
C\!\ell_{0,m}$ is an isometry for the standard scalar 
product on 
$\wedge {\Bbb R}^m = C\!\ell_m$. The following lemma is true by 
construction. 

\begin{lemma} The scalar 
product $<\cdot , \cdot >$ on $C\!\ell_{0,m}$ is invariant under 
left and right multiplications by unit vectors   $x \in {\Bbb R}^{0,m}$. 
In particular, 
if $m\equiv 0$ or $1\pmod{4}$, $h = {<\cdot ,\cdot >}$ is the 
(admissible) ${Pin} (0,2m)$-invariant  scalar product on the irreducible 
$C\!\ell_{0,2m}$-module $S = S_{0,2m} = C\!\ell_{0,m}$. 
\end{lemma} 

If $m \equiv 2$ or $3\pmod{4}$, we extend the scalar product 
${<\cdot ,\cdot >}$ on $C\!\ell_{0,m}$ to a symmetric complex bilinear form 
${<\cdot ,\cdot >}_{{\Bbb C}}$ on $S = \wedge {\Bbb C}^m$. 
Using the operator $c = c_{0,m}$ of complex conjugation with respect to 
the real form $C\!\ell_{0,m} = C\!\ell_m^0 + iC\!\ell_m^1$ of 
${\Bbb C}\!\ell_m$, we define a (real) scalar product 
$h = Re\, {<c\, \cdot ,\cdot >}_{{\Bbb C}}$ on $S$. 

\begin{lemma} Let $m \equiv 2$ or $3\pmod{4}$. Then  $h = 
Re\, {<c\, \cdot ,\cdot >}_{{\Bbb C}}$ is the (admissible) 
${Pin} (0,2m)$-invariant  scalar product on the irreducible 
$C\!\ell_{0,2m}$-module $S = {\Bbb C}\!\ell_m$. 
\end{lemma} 

Now we construct ($h$-)admissible bases  of  
the Schur algebra ${\cal C} = {\cal C}_{0,2m}$ for all the values of 
$m\pmod{4}$. 
\begin{prop} \label{02mSchurProp} If $m \equiv 0\pmod{4}$,  an admissible 
basis of the Schur 
algebra $${\cal C}_{0,2m} \cong {\Bbb R} \oplus {\Bbb R}$$ is 
given by the endomorphisms $ Id$ and $E=\alpha$ of $S = C\!\ell_{0,m}$: 
$\tau (E) = -1$, $\sigma (E) = \sigma_h(E) = +1$, $\iota (E) = +1$. 

If $m \equiv 1\pmod{4}$,  an admissible basis of 
${\cal C}_{0,2m} \cong {\Bbb C}$ is given by the endomorphisms $ Id$ and
$J = L_{\omega} \circ \alpha$ of $S = C\!\ell_{0,m}$ (where $\omega$ 
is a volume element of $C\!\ell_{0,m}$): $\tau (J) = -1$, $\sigma (J) 
= -1$. 

The space $\cal B$ of $\fr{so}(0,2m)$-invariant bilinear forms on $S$ is 
spanned by the admissible elements $h$ and $h_E$ if $m \equiv 0\pmod{4}$ and 
by $h$ and $h_J$ if $m \equiv 1\pmod{4}$. Their fundamental invariants 
$(\tau ,\sigma , \iota )$ are $(\tau ,\sigma , \iota )(h) = (+1,+1,+1)$, 
$(\tau ,\sigma , \iota )(h_E) = (-1,+1,+1)$ if $m \equiv 0\pmod{4}$ and 
$(\tau ,\sigma )(h) = (+1,+1)$, 
$(\tau ,\sigma )(h_J) = (-1,-1)$ if $m \equiv 1\pmod{4}$. 
\end{prop}

If $m \equiv 2$ or $3\pmod{4}$, we consider  the following operators on $S
= {\Bbb C}\!\ell_m$: 
\[ I: s\mapsto is\, ,\: J = L_{\omega}\circ c\, ,\: K= IJ \: \mbox{and} \: 
E =  \alpha \quad  (\omega = \omega_{0,m})\, .\] 
\begin{prop} Let $m \equiv 2$ or $3\pmod{4}$. The Schur algebra 
${\cal C}_{0,2m}$ ($\cong {\Bbb H}\oplus {\Bbb H}$ if 
$m \equiv 2\pmod{4}$ and $\cong {\Bbb C}(2)$ if $m \equiv 3\pmod{4}$) 
is generated by the admissible operators $I$, $J$ and $E$, which satisfy 
the following identities: 
\[ I^2 = J^2 = L_{\omega}^2 = -1\, ,\quad E^2 =c^2 = +1\, ,\] 
\[ \{ I,J\} = [I,E] = [I,L_{\omega}] = \{I,c\} = 0\, ,\] 
\[ [J, L_{\omega}] = [J,c] =[E,c] = [L_{\omega},c] = 0\, ,\] 
\[ [J,E] = [L_{\omega},E] = 0 \quad \mbox{if} \quad m \equiv 2\pmod{4}\, ,\] 
\[ \{ J,E\} = \{ L_{\omega},E\} = 0 \quad \mbox{if} 
\quad m \equiv 3\pmod{4}\, .\] 
An admissible basis of the Schur algebra is given by the endomorphisms 
$ Id$, $I$, $J$, $K$, $E$, $EI$, $EJ$, $EK$.  Their fundamental 
invariants $(\tau ,\sigma , \iota )$ are given in the next table, where 
the value of $m$ is modulo 4. \\ 

\begin{centering}
\begin{tabular}{|l||c|c|c|c|c|c|c|c|}\hline 
$m$: & $ Id$ & $I$ & $J$ & $K$ & $E$ & $EI$ & $EJ$ & $EK$\\\hline\hline
$2$ & $+++$ & $+-+$ & $--+$ & $--+$ & $-++$ & $--+$ & $+-+$ & $+-+$\\\hline 
$3$ & $+++$ & $+-+$ & $+--$ & $+--$ & $-++$ & $--+$ & $-+-$ & $-+-$\\\hline 
\end{tabular}
\vskip4pt
\end{centering}
The fundamental invariants  of the corresponding admissible basis for the 
space ${\cal B} = {\cal B}_{0,2m}$ (of $\fr{so}(0,2m)$-invariant 
bilinear forms on 
$S_{0,2m}$) are as follows: 
\begin{centering}
\vskip4pt
\begin{tabular}{|l||c|c|c|c|c|c|c|c|}\hline 
$m$: & $h$ & $h_I$ & $h_J$ & $h_K$ & $h_E$ & $h_{EI}$ & $h_{EJ}$ & 
$h_{EK}$\\\hline\hline
$2$ & $+++$ & $+-+$ & $--+$ & $--+$ & $-++$ & $--+$ & $+-+$ & $+-+$\\\hline 
$3$ & $+++$ & $+-+$ & $+--$ & $+--$ & $-++$ & $--+$ & $-+-$ & $-+-$\\\hline 
\end{tabular} 

\end{centering}
\end{prop}
  
\begin{thm} Every $\fr{so}(0,2m)$-equivariant  embedding 
${\Bbb R}^{0,2m} \hookrightarrow (S\otimes S)^{\ast}$,
$S = S_{0,2m}$, is a 
linear combination of the embeddings 
\[ j_{\rho}(h): {\Bbb R}^{0,2m}\hookrightarrow (S^+\vee S^-)^{\ast}
\quad \mbox{and} \quad j_{\rho}(h_E): {\Bbb R}^{0,2m}
\hookrightarrow (S^+\wedge S^-)^{\ast}\] 
if $m\equiv 0\pmod{4}$ and a 
linear combination of 
\[ j_{\rho}(h)\: \mbox{and} \: j_{\rho}(h_J): {\Bbb R}^{0,2m}
\hookrightarrow \vee^2S^{\ast} \quad \mbox{if} \quad m\equiv 1\pmod{4}\, .\] 

\noindent 
If $m\equiv 2$ or $3\pmod{4}$ every $\fr{so}(0,2m)$-equivariant  embedding 
${\Bbb R}^{0,2m}\hookrightarrow (S\otimes S)^{\ast}$ is a 
linear combination of the embeddings $j_A = j_{\rho}(h_A)$, $A\in {\cal C} = 
{\cal C}_{0,2m}$ 
admissible, whose image is contained in the dual of the subspaces indicated 
in Table \ref{02mTable} depending on $m\pmod{4}$. 
\begin{table}[ht]\caption[$\fr{so}(0,2m)$-equivariant  embeddings 
$j_A: {\Bbb R}^{0,2m} 
\hookrightarrow (S\otimes S)^{\ast}$]{\label{02mTable}
$\fr{so}(0,2m)$-equivariant  embeddings 
$j_A: {\Bbb R}^{0,2m} 
\hookrightarrow (S\otimes S)^{\ast}$} 

\begin{centering}
\begin{tabular}{|l||c|c|} \hline 
$j_{Id}$ & $S^+\vee S^-$ & $S^+\vee S^-$\\\hline 
$j_I$ & $S^+\wedge S^-$ & $S^+\wedge S^-$\\\hline 
$j_J$ & $S^+\vee S^-$ & $\wedge^2S^+ + \wedge^2S^-$\\\hline 
$j_K$ & $S^+\vee S^-$ & $\wedge^2S^+ + \wedge^2S^-$\\\hline 
$j_E$ & $S^+\wedge S^-$ & $S^+\wedge S^-$\\\hline 
$j_{EI}$ & $S^+\vee S^-$ & $S^+\vee S^-$\\\hline 
$j_{EJ}$ & $S^+\wedge S^-$ & $\wedge^2S^+ + \wedge^2S^-$\\\hline 
$j_{EK}$ & $S^+ \wedge S^-$ & $\wedge^2S^+ + \wedge^2S^-$\\\hline\hline 
$m$: & $2$ & $3$ \\\hline 
\end{tabular} 

\end{centering}
\vskip18pt
\end{table}  
\end{thm} 
\subsection{Case of odd dimension} 
Consider the 
orthogonal decomposition 
$$({\Bbb R}^{0,2m+1}, {<\cdot ,\cdot >}) = 
{\Bbb R}e_0 + {\Bbb R}^{0,2m},$$ where ${<e_0,e_0>} = -1$. 
Let $\rho$ 
denote the irreducible representation of $C\!\ell_{0,2m}$ on 
$S_{0,2m}$ defined in section \ref{negevenSubsec}. We  will extend
$\rho$ to an irreducible representation $\tilde{\rho}$ of $C\!\ell_{0,2m+1}$ 
on $S = S_{0,2m+1}$, where $S_{0,2m+1} = S_{0,2m}$ if $m\equiv 0$, $2$ 
or $3\pmod{4}$ and  $S_{0,2m+1} = S_{0,2m}\otimes {\Bbb C} = 
{\Bbb S}_{2m}$ if $m\equiv 1\pmod{4}$. If $m\equiv 2$ or 
$3\pmod{4}$, $S_{0,2m}= {\Bbb S}_{2m}$ admits the 
$C\!\ell_{0,2m}$-invariant complex structure $I$. For $m\equiv 1\pmod{4}$ 
multiplication by $i$ is a $C\!\ell_{0,2m}$-invariant complex structure
on $S_{0,2m}\otimes {\Bbb C}$ and will also be denoted by $I$. 

\begin{prop} The following formulas define an irreducible representation $\tilde{\rho}$ of $C\!\ell_{0,2m+1}$ 
on $S_{0,2m+1}$. 
\[ \tilde{\rho}|{\Bbb R}^{0,2m} = \rho |{\Bbb R}^{0,2m}\, ,\] 
\[ \tilde{\rho}(e_0) = \left\{ 
\begin{array}{r@{\quad \mbox{if} \quad m\equiv \:}l} 
\rho (\omega_{0,2m}) & 0\quad \mbox{or} \quad 2\pmod{4} \\
I \circ \rho (\omega_{0,2m}) & 1\quad \mbox{or} \quad 3\pmod{4}\, ,
\end{array} \right. \] 
where, in the case   $m\equiv 1\pmod{4}$,  $\rho$ has been extended 
complex linearly to a representation  on $S_{0,2m+1} = S_{0,2m}\otimes 
{\Bbb C}$.  $S =  S_{0,2m+1}$ is irreducible as 
$C\!\ell_{0,2m+1}^0$-module if $m\not\equiv 3\pmod{4}$ and the sum 
$S = S^+ + S^-$ of the two equivalent irreducible  
$C\!\ell_{0,2m+1}^0$-modules 
$S^+ = S^{\hat{J}}$ and $S^- = iS^{\hat{J}}$ if $m\equiv 3\pmod{4}$, 
where $S^{\hat{J}}$ is the fixed point set of a $\fr{so}(0,2m+1)$-invariant 
real structure $\hat{J}$ on $S$ (the explicit expression for $\hat{J}$ 
will be  given below). 
\end{prop} 

Next we describe the ${Pin} (0,2m+1)$-invariant scalar product 
$h = h_{0,2m+1}$ on $S =  S_{0,2m+1}$. Let $h_{0,2m}$ denote the 
${Pin} (0,2m)$-invariant 
scalar product  on $S_{0,2m+1}=S_{0,2m}$ if $m\equiv 0$, $2$ or 
$3\pmod{4}$ and by $h_{0,2m}^{\Bbb C}$ the complex bilinear 
extension of the  ${Pin} (0,2m)$-invariant 
scalar product  on $S_{0,2m}$ to a ${Pin} (0,2m)$-invariant complex 
bilinear form on $S_{0,2m+1} = {\Bbb S}_{2m} =  S_{0,2m}\otimes 
{\Bbb C}$ if $m\equiv 1\pmod{4}$. 
\begin{lemma} The ${Pin} (0,2m+1)$-invariant scalar product $h = h_{0,2m+1}$ 
on $S =  S_{0,2m+1}$ is given by $h = h_{0,2m}$ if $m\equiv 0$, $2$ or 
$3\pmod{4}$ and by $h = Re\, h_{0,2m}^{{\Bbb C}}(c\, \cdot ,\cdot )$ if 
$m\equiv 1\pmod{4}$, where $c$ is complex conjugation with respect to 
$S_{0,2m}\subset S_{0,2m}\otimes {\Bbb C}$.
\end{lemma} 

If $m\not\equiv 0\pmod{4}$, we have on $S_{0,2m+1} = {\Bbb C}\!\ell_m = 
C\!\ell_{0,m} + iC\!\ell_{0,m}$ the operator $c = c_{0,m}$ of complex 
conjugation. 
Using it we define an endomorphism $\hat{J}$ of $S_{0,2m+1} = 
{\Bbb C}\!\ell_m$ by 
\[ \hat{J} := L_{\omega}\circ \alpha \circ c\, ,\] 
where $\omega = \omega_{0,m}$ is a volume element of $C\!\ell_{0,m}$ and 
$\alpha | {\Bbb C}\!\ell_m^0 = +{Id}$, $\alpha | 
{\Bbb C}\!\ell_m^1 = -{Id}$.
\begin{prop} Let $m\not\equiv 0\pmod{4}$. The Schur 
algebra ${\cal C} = {\cal C}_{0,2m+1}$ is generated by the endomorphisms 
$I$ and $\hat{J}$ of $S = S_{0,2m+1} = 
{\Bbb C}\!\ell_m$, which satisfy the following relations: 
$I^2 = -1$, $\{ I,\hat{J} \} = 0$. Moreover, $\hat{J}^2 = +{Id}$ and 
${\cal C}_{0,2m+1}\cong {\Bbb R}(2)$ if $m\equiv 3\pmod{4}$ and 
$\hat{J}^2 = -{Id}$ and ${\cal C}_{0,2m+1}\cong {\Bbb H}$ if 
$m\equiv 1$ or $2\pmod{4}$. An admissible basis of ${\cal C}_{0,2m+1}$ 
is given by the endomorphisms $ Id$, $I$, $\hat{J}$ and $\hat{K} = I\hat{J}$. 
Their fundamental invariants $(\tau ,\sigma , \iota )$ together with 
the invariants of the associated admissible basis for  the space $\cal B$ 
of $\fr{so}(0,2m+1)$-invariant bilinear forms are given in the  Table 
\ref{02m+1fiTable} ($\iota$ is only defined if $m\equiv 3\pmod{4}$).
\begin{table}[ht]\caption[Fundamental invariants of admissible endomorphisms 
and bilinear forms of $S_{0,2m+1}$]{\label{02m+1fiTable}Fundamental invariants of admissible endomorphisms and bilinear forms of $S_{0,2m+1}$}  

\begin{centering}
\begin{tabular}{|l||c|c|c|c||c|c|c|c|}\hline 
$m$: & $ Id$ & $I$ & $\hat{J}$ & $\hat{K}$ & 
$h$ & $h_I$ & $h_{\hat{J}}$ & $h_{\hat{K}}$ \\\hline\hline
$1$ & $++$ & $+-$ & $--$ & $--$ & 
$++$ & $+-$ & $--$ & $--$\\\hline  
$2$ & $++$ & $+-$ & $+-$ & $+-$ & 
$++$ & $+-$ & $+-$ & $+-$\\\hline  
$3$ & $+++$ & $+--$ & $-++$ & $-+-$ & 
$+++$ & $+--$ & $-++$ & $-+-$\\\hline 
\end{tabular}

\end{centering}
\vskip18pt
\end{table} 
If $m\equiv 0\pmod{4}$,   ${\cal C}_{0,2m+1} = {\Bbb R} {Id}$. 
\end{prop} 

\begin{thm} Every $\fr{so}(0,2m+1)$-equivariant  embedding 
${\Bbb R}^{0,2m+1}\hookrightarrow (S\otimes S)^{\ast}$ 
is proportional to $j_{\rho}(h): {\Bbb R}^{0,2m+}\hookrightarrow 
\vee^2S^{\ast}$ if $m\equiv 0\pmod{4}$ and a linear combination of the 
embeddings $j_A = j_{\rho}(h_A)$, $A = Id$, $I$, $\hat{J}$ and $\hat{K}$ 
if $m\not\equiv 0\pmod{4}$. The image of the $j_A$ is contained in the dual 
of the subspaces indicated in Table \ref{02+1mTable}
\begin{table}[ht]\caption[$\fr{so}(0,2m+1)$-equivariant  embeddings 
$j_A: {\Bbb R}^{0,2m+1} 
\hookrightarrow (S\otimes S)^{\ast}$]{\label{02+1mTable}
$\fr{so}(0,2m+1)$-equivariant  embeddings 
$j_A: {\Bbb R}^{0,2m+1} 
\hookrightarrow (S\otimes S)^{\ast}$} 

\begin{centering}
\begin{tabular}{|l||c|c|c|} \hline
$j_{Id}$ & $\vee^2S$ & $\vee^2S$ & $S^+\vee S^-$\\\hline
$j_J$ &  $\wedge^2S$ &  $\wedge^2S$ & $S^+\wedge S^-$\\\hline 
$j_{\hat{J}}$ & $\vee^2S$ & $\wedge^2S$ & $\wedge^2S^+ +\wedge^2S^-$\\\hline 
$j_{\hat{K}}$ & $\vee^2S$ & $\wedge^2S$ & $\wedge^2S^+ +\wedge^2S^-$
\\\hline\hline 
$m$: & $1$ & $2$ & $3$ \\\hline 
\end{tabular} 

\end{centering}
\vskip18pt
\end{table}  
\end{thm}  
\section{Complete classification} 
Every pseudo Euclidean space $V$ admits a (unique up to an isometry) 
 orthogonal decomposition 
$V = V_1 + V_2$, where $V_1 = {\Bbb R}^{m,m}$ and the scalar 
product of $V_2$ is 
positively or negatively defined. Now we consider the case when $V_1 \neq 0$ 
and $V_2 \neq 0$, the other cases were treated in the sections 
\ref{signmmsubSec}, \ref{posSec} and \ref{negSec}.   
We denote by $S_i$, $i=1,2$, the  irreducible $C\!\ell (V_i)$-module 
constructed in the sections \ref{signmmsubSec} and  \ref{posSec}, 
\ref{negSec} respectively. Then $S = S_1 \otimes S_2$ carries the structure 
of irreducible module for the Clifford algebra $C\!\ell (V) = 
C\!\ell (V_1) \hat{\otimes} C\!\ell (V_2)$, s.\ Proposition \ref{redProp}. 
By  Proposition \ref{b1b2Prop}, to every admissible bilinear form 
$\beta_2$ (respectively endomorphism $A_2$) on $S_2$ we associate an 
admissible bilinear form  $\beta = \beta_1 \otimes \beta_2$ (respectively 
endomorphism $A_1\otimes A_2$) on $S$.  In the sections 
\ref{posSec} and \ref{negSec} we have contructed admissible bases for the 
space ${\cal B}_2$ of $\fr{so}(V_2)$-invariant bilinear forms on $S_2$ and 
for the Schur algebra ${\cal C}_2$ of $S_2$. Therefore, this explicit 
correspondence defines an injective linear mapping 
$\phi : \beta_2 \mapsto \beta = 
\phi (\beta_2)$ (respectively $\psi: A_2 \mapsto A = \psi (A_2)$) from 
${\cal B}_2$ into the space $\cal B$ of $\fr{so}(V)$-invariant bilinear 
forms on $S$ (respectively from ${\cal C}_2$ into the Schur algebra 
$\cal C$ of $S$).  Moreover, $\phi$ and $\psi$ are actually isomorphisms, 
because the Schur algebras of $S$ and $S_2$ are isomorphic, due to the fact
that $V$ and $V_2$ have the same signature $s$, s.\ Corollary \ref{SchurCor}. 
So we have essentially proved: 
\begin{thm}\label{basicThm} There exist natural   isomorphisms $\phi : {\cal B}_2 
\rightarrow  {\cal B}$ of vector spaces and $\psi : {\cal C}_2 
\rightarrow  {\cal C}$ of algebras mapping admissible elements onto admissible
elements. Under these maps, the fundamental invariants of admissible elements 
transform according to the rules given in Proposition \ref{multiplProp}. 

In particular, if $m\equiv 0\pmod{4}$, then   $\phi$ and $\psi$ preserve the 
fundamental invariants ({\bf (4,4)-periodicity}). 
\end{thm} 

\noindent 
{\bf Proof:} We recall that by Proposition \ref{ffeProp} the Schur algebra
${\cal C}_{m,m}$ of $S_1 = S_{m,m}$ has the admissible basis $({Id}, E)$ 
and $E^2 = +{Id}$. This implies    that the  vector space isomorphism 
$\psi$ is actually an isomorphism of {\em algebras}. The 
(4,4)-periodicity follows from 
\[ \sigma (f_E) = \iota (f_E) = \sigma_f(E) = \sigma_{f_E}(E) = \iota (E) 
= +1\, .\quad \Box\] 
Recall that ${\cal B}_{p,q}$ denotes  the space of $\fr{so}(p,q)$-invariant 
bilinear forms on the $\fr{so}(p,q)$ spinor module $S_{p,q}$ and 
${\cal C}_{p,q}$ is the Schur algebra of $S_{p,q}$.  
\begin{cor} {\bf ((8,0)- and (0,8)-periodicity)} There exist natural 
isomorphisms 
\[ \phi_{8,0} : {\cal B}_{p,q} \rightarrow {\cal B}_{p+8,q}\quad 
\mbox{and} \quad \phi_{0,8} : {\cal B}_{p,q} \rightarrow {\cal B}_{p,q+8}\]
of vector spaces and 
\[ \psi_{8,0} : {\cal C}_{p,q} \rightarrow {\cal C}_{p+8,q}\quad 
\mbox{and} \quad \psi_{0,8} : {\cal C}_{p,q} \rightarrow {\cal C}_{p,q+8}\]
of algebras mapping the admissible elements onto admissible elements  
preserving their fundamental invariants. 
\end{cor} 

\noindent 
{\bf Proof:} By Theorem \ref{basicThm}  ${\cal B}_{p,q}$ and ${\cal C}_{p,q}$
have admissible bases. Now we recall from sections \ref{posSec} and 
\ref{negSec} that if $k\equiv 0\pmod{8}$, then ${\cal C}_{k,0} \cong 
{\cal C}_{0,k}$ has an admissible basis, which was denoted by  $({Id}, E)$, 
such that $(\tau ,\sigma , \iota )(E) = (-1,+1,+1)$ and, of course, 
$(\tau ,\sigma , \iota )(Id) = (+1,+1,+1)$. The existence of the maps 
$\psi_{8,0}$ and  $\psi_{0,8}$ follows from $\tau (Id) \iota (Id) = -
\tau (E) \iota (E)$. They preserve the fundamental invariants, because 
$\sigma (Id) =  \iota (Id) = \sigma (E) = \iota (E) = +1$. The existence and 
properties of $\phi_{8,0}$ and  $\phi_{0,8}$ are proved similarly. $\Box$ 

\begin{cor} Every $\fr{so}(V)$-equivariant mapping $j: 
V\rightarrow (S\otimes S)^{\ast}$ is a linear combination of the 
embeddings $j_A = j_{\rho}(h_A)$, where $h$ is the canonical bilinear form 
on the spinor module $S$ of $\fr{so}(V)$ and $A$ are admissible elements of
the Schur algebra $\cal C$ of $S$. 
\end{cor} 

To obtain an overview over all possible $N$-extended Poincar\'e algebras 
$\fr{p}(V) + S$, $N = \pm 1, \pm 2$, it is useful to define 
the invariants $\sigma$ and $\iota$ for embeddings $j: 
V\hookrightarrow (S\otimes S)^{\ast}$ having special properties. 
More precisely, we put $\sigma (j) = +1$ if $jV\subset \vee^2S^{\ast}$ 
and  $\sigma (j) = -1$ if $jV\subset \wedge^2S^{\ast}$. If $S = S^+ + S^-$, 
we  define $\iota (j) = +1$ if 
$jV\subset (S^+\otimes S^+ + S^-\otimes S^-)^{\ast}$ and 
$\iota (j) = -1$ if 
$jV\subset (S^+\otimes S^-)^{\ast}$. 

Note that the fundamental invariants of $j_A = j_{\rho}(h_A)$, 
$A\in {\cal C}$ admissible,
are easily computable: 
\[ \sigma (j_A) = \tau (h_A) \sigma (h_A) = \tau (h) \tau (A) \sigma (h) 
\sigma (A) \quad \mbox{and} \quad 
\iota (j_A) = -\iota (h_A) = -\iota (h) \iota (A)\, .\] 

Recall that $\cal J$ denotes the space of $\fr{so}(V)$-equivariant 
mappings $j: V\rightarrow (S\otimes S)^{\ast}$. We define the subspaces
\[ {\cal J}^{\sigma_0} := \{ j\in {\cal J}| \sigma (j) = 
\sigma_0\} \cup \{ 0\} \quad \mbox{and} \]
\[ {\cal J}^{\sigma_0\iota_0}:= \{ j\in {\cal J}^{\sigma_0}| \iota (j) = 
\iota_0\} \cup \{ 0\} \] 
and put 
\[ L^{\sigma_0}:= \dim {\cal J}^{\sigma_0}\, , \quad   
 L^{\sigma_0\iota_0}:= \dim {\cal J}^{\sigma_0\iota_0}\, .\] 
We shall write $L^+$, $L^{+-}$, ... instead of the more cumbersome 
$L^{+1}$, $L^{+1\, -1}$, ... 

Remark that $L^+$ ($= L^{++} + L^{+-}$ if 
$S = S^+ + S^-$) is the maximal number of linearly independent 
super algebra structures on $\fr{p}(V) + S$ and that 
$L^-$ ($= L^{-+}+ L^{--}$) is the number of $\Bbb Z_2$-graded 
Lie algebra structures on $\fr{p}(V) + S$. 

\begin{thm} The numbers $(L^+,L^-)$ and $(L^{++},L^{+-},L^{-+},L^{--})$ 
depend only on the dimension $n = \dim V = p+q$ and the signature 
$s = p-q$ of $V = {\Bbb R}^{p,q}$ modulo 8. Moreover, 
they admit the  {\bf mirror super symmetry} $n\mapsto -n$. 
More precisely, 
\begin{eqnarray*}
 L^+(-n,s) &=& L^-(n,s)\quad \mbox{and}\\ 
 L^{+\, \iota_0} (-n,s) &=& L^{- \, \iota_0}(n,s)\, , \quad \iota_0 = \pm \, .
\end{eqnarray*}
Their values are given in Table \ref{mirrorsymmTable}.     
\begin{table}[ht]\caption[Numbers of different types of extended Poincar\'e 
algebras $\fr{p}(p,q) + S_{p,q}$ depending on $n= p+q$ and 
$s = p-q$ modulo 8]{\label{mirrorsymmTable}Numbers of  extended Poincar\'e 
algebras $\fr{p}(p,q) + S_{p,q}$ of different types  depending on $n= p+q$ and 
$s = p-q$ modulo 8}     

\begin{centering}
\begin{tabular}{|r||c|c|c|c|c|c|c|c|}\hline 
$s$: & \multicolumn{8}{c|}{$(L^{++},L^{+-},L^{-+},L^{--})(n,s)$ or 
$(L^+,L^-)(n,s)$}\\ \hline\hline 
4 & & 2,0,6,0 & & 0,4,0,4 & & 6,0,2,0 & & 0,4,0,4\\ \hline
3&1,3& &1,3& &3,1& & 3,1& \\ \hline
2& & 0,2,4,2 & & 2,2,2,2 & &4,2,0,2& &2,2,2,2\\ \hline
1 &0,1,2,1& &0,1,2,1& &2,1,0,1& &2,1,0,1& \\ \hline
0& &0,0,2,0& &0,1,0,1&  &2,0,0,0& &0,1,0,1\\ \hline
-1 &0,1& & 0,1& &1,0& &1,0& \\ \hline
-2& &0,2& &1,1&  &2,0& &1,1\\ \hline
-3 &1,3& &1,3& &3,1& &3,1& \\ \hline\hline 
$n$: &-3&-2&-1&0&1&2&3&4\\ \hline   
\end{tabular}    

\end{centering}
\vskip18pt
\end{table}  
\end{thm} 

\noindent 
{\bf Proof:} This follows from Theorem \ref{basicThm} and 
the tables of  sections \ref{signmmsubSec}, \ref{posSec} and \ref{negSec} 
by straightforward computation. $\Box$

In the complex case we consider the space ${\cal J}_c$ of 
$\fr{so}(m,{\Bbb C})$-equivariant mappings 
${\Bbb C}^m \rightarrow   
({\Bbb S}_m\otimes {\Bbb S}_m)^{\ast}$ and define the 
invariants $\sigma$, $\iota$ and  
the spaces ${\cal J}^+_c$,  ${\cal J}^{+-}_c$ etc.\  
as in the real case ($\iota$ is only defined if the complex 
$\fr{so}(m,{\Bbb C})$ spinor 
module ${\Bbb S}_m$ is reducible ${\Bbb S}_m = {\Bbb S}_m^+ + 
{\Bbb S}_m^-$). Their dimensions are  denoted by $L^+_c$, $L^{+-}_c$ 
etc.\  

\begin{thm} The numbers $(L_c^+,L_c^-)$ and $(L_c^{++},L_c^{+-},L_c^{-+},L_c^{--})$ depend only on $m\pmod{8}$. Moreover, 
they admit the mirror super symmetry $m\mapsto -m$. 
More precisely, 
\begin{eqnarray*}
 L_c^+(-m) &=& L_c^-(m)\quad \mbox{and}\\ 
 L_c^{+\, \iota_0} (-m) &=& L_c^{- \, \iota_0}(m)\, , \quad \iota_0 = \pm \, .
\end{eqnarray*}
Their values are given in the next table. \\ 

\begin{centering}
\begin{tabular}{|r||c|c|c|c|c|c|c|c|} \hline 
   & $0,1$ & $0,0,2,0$ & $0,1$ & $0,1,0,1$ & $1,0$ & $2,0,0,0$ & $1,0$ & 
$0,1,0,1$  
\\\hline\hline 
$m$: & $-3$ & $-2$ & $-1$ & $0$ & $1$ & $2$ & $3$ & $4$ \\\hline  
\end{tabular} 

\end{centering}
\end{thm} 

\noindent
{\bf Proof:} follows from section  \ref{cxcaseSec}. $\Box$  

\end{document}